\documentclass[a4paper,12pt]{article}
\usepackage{amssymb}
\usepackage{graphicx}
\newcommand{\M}{{\cal M}}

\newcommand{\RR}{ I\!\!R}

\newcommand{\R}{\mathbb{R}}

\newcommand{\N}{\mathbb{N}}

\newcommand{\beq}{\begin{equation} }
\newcommand{\eqq}{\end{equation} }
\newcommand{\cuad}{{\sqcap\kern-.68em\sqcup}}

\newcommand{\norm}[1]{\|#1\|}
\newcommand{\equ}[1]{(\ref{#1})}
\newcommand\rn{\mathbb{R}^N}

\newtheorem{teo}{Theorem}[section]
\newtheorem{prop}{Proposition}[section]
\newtheorem{proposition}{Proposition}[section]
\newtheorem{lema}{Lemma}[section]
\newtheorem{lemma}{Lemma}[section]

\newtheorem{remark}{Remark}[section]
\newcommand{\bremark}{\begin{remark} \em}
\newcommand{\eremark}{\end{remark} }

\def\linf{L^\infty(\Omega)}
\def\lp{L^p(\Omega)}

\def\mm{{{\cal M}^-_{\lambda,\Lambda}}}

\newcommand{\lpl}{\lambda_1^+}
\newcommand{\lmin}{\lambda_1^-}

\newcommand{\lplf}{\lambda_1^+(F, \Omega)}
\newcommand{\lminf}{\lambda_1^-(F,\Omega)}

\newcommand{\fipl}{\varphi_1^+}
\newcommand{\fimin}{\varphi_1^-}

\def\beeq{\begin{equation}}
\def\eeq{\end{equation}}
\newcommand{\begeqaet}{\begin{eqnarray*}}
\newcommand{\eneqaet}{\end{eqnarray*}}
\newcommand{\calc}{{\cal C}}

\newcommand{\s}{{\tau}}

\newcommand{\re}[1]{(\ref{#1})}
\begin{document}
\begin{center}{\bf\Large  Resonance phenomena for second-order stochastic control equations }\medskip

{Patricio FELMER}, {Alexander QUAAS}, {Boyan SIRAKOV}\end{center}




\begin{abstract}
We study the existence and the properties of solutions to the Dirichlet problem for uniformly elliptic second-order Hamilton-Jacobi-Bellman operators, depending on  the principal eigenvalues of the operator.
\end{abstract}
\date{}

\setcounter{equation}{0}
\section{ Introduction}

In this article we consider the Dirichlet problem
\beq\label{first}
\left\{\begin{array}{rclcc} F(D^2u,Du,u,x)&=&f(x)&\mbox{in}& \Omega,\\
u&=&0&\mbox{on}& \Omega, \end{array}\right.\eqq where the  second
order differential operator $F$ is of Hamilton-Jacobi-Bellman (HJB)
type and $\Omega\subset\rn$ is a bounded domain. These equations --
see the book \cite{FS} and the surveys \cite{K1}, \cite{Son},
\cite{Ca}, as well as \cite{L} (various other references will be given below) -- have
been very widely studied because of their connection with the
general problem of Optimal Control for Stochastic Differential
Equations (SDE). We recall that a powerful approach to  this problem
is the so-called Dynamic Programming Method, originally due to R.
Bellman, which indicates that the optimal cost (value) function  of
a controlled SDE should be a solution of a PDE like \re{first}. More
precisely, let us have a stochastic process $X_t$ satisfying
$$
dX_t = b^{\alpha_t}(X_t) dt + \sigma^{\alpha_t}(X_t)dW_t\,,
$$
with $X_0 = x$, for some $x\in\Omega$, and the cost function
$$
J(x,\alpha) = \mathbb{E} \int_0^{\tau_x} f(X_t) \exp \{\int_0^t c^{\alpha_s}(X_s) ds\}\,dt,
$$
where $\tau_x$ is the first exit time from ${\Omega}$ of $X_t$, and $\alpha_t$ is an index (control) process with values in a set ${\cal A}$. Then the optimal cost function
$v(x) = \inf_{\alpha\in {\cal A}} J(x,\alpha)$ is such that $-v$ is  a solution of \re{first}, which is in the form
\beq\label{equ1}
\left\{\begin{array}{rclcc}
\displaystyle\sup_{\alpha\in {\cal A}} \{\mathrm{tr} (A^\alpha(x)D^2u) + b^\alpha(x).Du
+ c^\alpha(x) u \} &=& f(x) &\mbox{ in }& \Omega\\ u&=&0 &\mbox{on}& \partial
\Omega.\end{array}\right.
 \eeq
We are going to study this boundary value problem under the following
hypotheses, which will be kept throughout the paper :  for some constants $0<\lambda\le \Lambda$, $\gamma\ge0$, $\delta\ge0$,
we assume $A^\alpha(x) := \sigma^\alpha(x)^T\sigma^\alpha(x)\in
C(\overline{\Omega})$, $\lambda I\! \le\! A^\alpha(x) \le\! \Lambda I$,
$|b^\alpha(x)|\le \gamma$,  $|c^\alpha(x)|\le \delta$, for almost all $
x\in \Omega$ and all $\alpha\in {\cal A}$, and $f\in L^p(\Omega)$, for some $p> N$.
We stress however that all our results are new even for operators with smooth coefficients.

Our main statements on resonance, applied to this setting, imply
in particular  that for some $A,b,c$ the optimal cost becomes arbitrarily
large or small, depending on the function $f$ which stays bounded. We give conditions under which \re{equ1} is solvable or not, and describe  properties of its solutions.

The majority of works on HJB equations concern {\it proper} equations,
that is, cases when $F$ is monotone in the variable $u$ ($c^\alpha\le 0$),
in which no resonance phenomena can arise. It was shown in the well-known papers
\cite{E1}, \cite{E2} and \cite{L2} that a proper equation of type \re{equ1}
has a unique strong solution, which is classical, if the coefficients are
smooth. Uniqueness in the viscosity sense was proved in \cite{J}, \cite{CIL},
\cite{CCKS} and \cite{Sw}.

Two of the authors recently showed in \cite{QS2} that  existence and uniqueness
 of viscosity solutions holds for a larger class of operators, including
nonproper operators whose principal eigenvalues -- defined below -- are
positive. This had been  proved much earlier for HJB operators with smooth
coefficients in \cite{L}, through a mix of probabilistic and analytic
techniques. Very recently existence, non-existence and multiplicity
results for cases when the eigenvalues are negative or have different signs,
 but are different from zero, appeared in \cite{A} and \cite{S1}.

 Thus the only situations which remain completely unstudied are the cases when  \re{equ1} is "at resonance", that is, when one of the principal eigenvalues of $F$ is zero. The present paper is devoted to this problem. We also obtain a number of new results for cases without resonance.

We shall make essential use of the work \cite{QS2}, where the properties of the eigenvalues are studied. In particular, based on the definition for the linear case in \cite{BNV}, it is shown in \cite{QS2} that the numbers
$$
\lambda_1^+(F,\Omega)=\sup\{\lambda\,|\, \Psi^+(F,\Omega,\lambda)\not =\emptyset\},\quad
\lambda_1^-(F,\Omega)=\sup\{\lambda\,|\, \Psi^-(F,\Omega,\lambda)\not =\emptyset\},
$$
where the sets $\Psi^+(F,\Omega,\lambda)$ and $ \Psi^-(F,\Omega,\lambda)$ are
defined  as
$$
\Psi^\pm(F,\Omega,\lambda)=\{\psi\in C(\overline{\Omega})\,|\,
\pm(F(D^2\psi,
D\psi,\psi,x)+\lambda\psi)\le 0,\;  \pm\psi>0
\mbox{ in } \Omega\},
$$
 are simple and isolated eigenvalues of $F$,  associated to a positive
and a negative eigenfunctions
$\varphi_1^+$, $\varphi_1^-\in W^{2,q}(\Omega)$, $q<\infty$,  and that their positivity
guarantees the validity of one-sided Alexandrov-Bakelman-Pucci type estimates -- see the review in
the next section.  From the optimal control point of view, $\lambda_1^+$
can be seen as the fastest exponential rate at which paths can exit the
domain, and $\lambda_1^-$ is the slowest one, we refer to the exact formulae given
in equalities (30)-(31) of \cite{L}. For extensions and related results on eigenvalues for fully nonlinear operators we refer  to  \cite{IY}, \cite{A}, where Isaacs operators are studied,  and to \cite{BD1}, \cite{BD2}, where more general singular fully nonlinear elliptic operators  are considered. When no confusion arises, we
write $\lambda_1^\pm$ or $\lambda_1^\pm(F)$, and we  always suppose
that $\lambda_1^+<\lambda_1^-$
 --- note it easily follows from the results in \cite{QS2} that $\lpl= \lmin$ can only happen if all  linear operators which appear in \re{equ1} have the same principal eigenvalues {\it and} eigenfunctions. For simplicity, we assume
that $\Omega$ is regular, in the sense that it  satisfies an uniform
interior ball condition, even though many of the results can be extended to general bounded domains.

We make the convention that all (in)equalities in the paper are meant to
hold in  the $L^p$-viscosity sense, as defined and studied in \cite{CCKS}.
Note however that it is known that any viscosity solution of \re{equ1} is
in $W^{2,p}(\Omega)$ and that any $W^{2,p}$-function which satisfies
\re{equ1} in the viscosity sense is also a strong solution,  that is,
it satisfies \re{equ1} a.e. in $\Omega$, see \cite{C}, \cite{CCKS},
\cite{Sw}, \cite{W}. All constants in the estimates will be allowed to depend on $N, \lambda,\Lambda,\gamma, \delta$ and $\Omega$.

Given a fixed function $h\in L^p(\Omega)$ which is not a multiple of the principal eigenfunction $\varphi_1^+$, along the paper we write
\begin{equation}\label{f}
f=t\varphi_1^++h, \quad t\in\R,
\end{equation}
and consider $t$ as a parameter. We note that all results and proofs below
hold without modifications if the function $\varphi_1^+$ in
 $\equ{f}$ is replaced by any other positive function,
 which vanishes on $\partial\Omega$ and whose interior normal derivative on the boundary is strictly positive.
We visualize  the set $S$ of solutions of \equ{equ1} in the
space
$C(\overline \Omega)\times \R$ as follows: $(u,t)\in {S}$ if
and only if $u$ is a solution of \equ{equ1} with $f=t\varphi_1^++h$.
The following notation will be useful:
given a subset $A\subset C(\overline \Omega)\times \R$ and $t\in\R$
we define $A_t=\{u\in C(\overline \Omega)\,|\, (u,t)\in A\}$ and $A_I=
\cup_{t\in I}A_t$, if $I$ is an interval.

Our purpose is to describe the set  ${S}=\{(u_t,t)\,|\, t\in\R \}$. When $\lambda_1^+(F)>0$
this can be done in a rather precise way, thanks to the results in
\cite{L},  \cite{QS2}.
\begin{teo}
\label{teo1}
Assume $\lambda_1^+(F)>0$. Then

\begin{itemize}
\item[1.] (\cite{QS2}) For every $t\in\R$ equation \equ{equ1} possesses
exactly one
solution $u=u_t$. In addition, if $f=t\varphi_1^++h\not =0$ and $f\le(\ge) 0$ then $u>(<)0$ in $\Omega$.
If $t<s$ then $u_t>u_s$ in $\Omega$.

\item[2.] The set $S$ is a Lipschitz continuous curve such that $t\to u_t(x)$ is convex for each $x\in\Omega$. There exist numbers $t^\pm=t^\pm(h)$ such that if $t\ge t^+$ ($t\le t^-$) then $u_t<(>)0$ in $\Omega$.  Moreover, for each compact $K\subset\subset \Omega$
$$
\lim_{t\to -\infty}\min_{x\in K}u_t(x)=+\infty\qquad\mbox{and} \qquad
\lim_{t\to +\infty}\max_{x\in K}u_t(x)=-\infty.
$$
\end{itemize}
\end{teo}
Next, we state our first main theorem, which   describes the set $S$ when the first eigenvalue is zero.
 In this case the set of solutions is  again a unique continuous curve,
but it exists only on a half-line with respect to $t$, and becomes
unbounded when $t$ is close or equal to a critical number $t^*_+$ -- see Figure 1 below.
Note the picture is very different from the one we obtain in the linear case - if $L$ is a linear operator then the Fredholm alternative for $Lu + \lambda_1(L) u = t\varphi_1(L) + h$ says this equation has a solution only for one value of $t$, and then any two solutions differ by a multiple of $\varphi_1(L)$.
\begin{teo}
\label{teo2}
Assume
$\lambda_1^+(F)=0$. Then

\begin{itemize}
\item[1.] There exists a number $t^*_+=t^*_+(h)$ such that if $t<t^*_+$ then there is no solution of \equ{equ1}, while for $t>t^*_+$ \equ{equ1} has a solution.

\item[2.] The set $S$ is a continuous curve  such that
$S_t$ is a singleton  for all
$t>t^*_+$, that is, solutions are unique for $t>t^*_+$. If $t^*_+\le t<s$ and $(u_t,t), (u_s,s)\in S$
then $u_t>u_s$ in $\Omega$. The map $t\to u_t(x)$ is convex for each $x\in\Omega$.

\item[3.] There exists $t^+=t^+(h)> t^*_+$ such that if $t\ge t^+$ then $u_t<0$
in $\Omega$, and for every compact $K \subset\subset\Omega$ we have $\displaystyle
\lim_{t\to +\infty}\max_{x\in K}u_t(x)=-\infty.
$

\item[4.] If $t=t^*_+$ then either:

~~~~~~~(i) Equation \equ{equ1} does not have a solution, that is,
$S_{t^*}$ is empty, $\lim_{t\searrow t^*_+}\min_{x\in K}u_t=+\infty$ for every compact $K \subset\subset\Omega$, and there exists $\epsilon=\epsilon(h)>0$ such that if $t\in (t^*_+,t^*_++\epsilon)$ then $u_t>0$, or

~~~~~~(ii) There exists a function $u^*$ such that $S_{t^*_+}= \{u^*+s\varphi_1^+\;|\; s\ge0\}$.
\end{itemize}
\end{teo}

In case the  two eigenvalues have opposite signs, a multiplicity phenomenon occurs. This situation was studied in  \cite{S1} and we recall it here.

\begin{teo}
\label{teo3}(\cite{S1})
Assume $\lambda_1^+(F)<0<\lambda_1^-(F)$. Then
there exists a number $t^*=t^*(h)$ such that
\begin{itemize}
\item[1.] If $t<t^*$ then there is no solution of \equ{equ1} ;
\item[2.] If $t>t^*$ then there are at least two solutions of
\equ{equ1} ; more precisely, for $t\in (t^*,\infty)$ there is a continuous curve of minimal solutions $u_t$ of \re{equ1} such that $t\to u_t(x)$ is convex and strictly decreasing for  $x\in\Omega$, and a  connected set of solutions different from the minimal ones ;
\item[3.]  If $t=t^*$ then there is at least one solution of \equ{equ1}.
\end{itemize}
\end{teo}

Note that in \cite{S1} the  properties of the two branches were
not described~; however by using the results there, our Lemma \ref{convexx} and
some topological and degree arguments, like in Sections \ref{sect3}-\ref{sect5}, they can be  obtained easily.

\medskip

Now we state our second main theorem, which describes properties
of the  set $S$ when the second eigenvalue is
 at resonance, that is, when  $\lambda_1^-(F)=0$. Here the analysis is more difficult
than in Theorem 1.2,
but still  the picture is quite clear.
\begin{teo}
\label{teo4}
Assume $\lambda_1^-(F)=0$. Then there exists $t^*_-=t^*_-(h)$ such that
\begin{itemize}
\item[1.] If $t<t^*_-$ then there is no solution of \equ{equ1}.
\item[2.] There is a closed connected set
$\calc\subset S$, such that
$\calc_t\not =\emptyset$ for all
$t>t^*_-$.

\item[3.] The set $S_I$ is bounded in $W^{2,p}(\Omega)$, for each compact
$I\subset (t^*_-,\infty)$.

\item[4.]  If we denote  $\alpha_t=\inf\{\sup_\Omega u\,|\, u\in S_t\}$, we have
$
\lim_{t\to +\infty}\alpha_t= +\infty.
$
\item[5.]  The set
$\calc_{[t^*_-, t^*_-+\varepsilon)}$
is unbounded in $\linf$, for all $\varepsilon>0$;  there exists $C=C(h)>0$
such that if $u\in S_{[t^*_-, t^*_-+\varepsilon)}$ and
$\norm{u}_{\linf}\ge C$ then $u<0$ in $\Omega$; if $u_n \in S_{[t^*_-, t^*_-+\varepsilon)}$ and $\norm{u_n}_{\linf}\to\infty$ then   $\max_K u_n \to -\infty$ for each compact $K\subset\Omega$.

\item[6.] If $S_{t^*_-}$ is unbounded in $\linf$ then there exists a function $u_*$ such that $S_{t^*_-}= \{u_*+s\varphi_1^-\;|\; s\ge0\}$.
\end{itemize}
\end{teo}

Both Theorems \ref{teo2} and \ref{teo4} are proved by a careful analysis of the behaviour of the sets of solutions to equations with positive (resp. negative) eigenvalues, when $\lpl(F)\searrow 0$ (resp. $\lmin(F)\nearrow 0$).

We note that not much is known on solutions of \re{equ1}  when
both eigenvalues are negative. Thus, before proving Theorem \ref{teo4}, we need
to analyze solutions of problems in which $\lmin(F)$ is small and negative. This is the content of the next theorem, which is of clear independent interest.

\begin{teo}\label{teo5}
There exists $0<L\le\infty$ such that if $\lambda_1^-(F)\in (-L,0)$ then
\begin{itemize}
\item[1.] There exists a
closed connected set
$\calc\subset S$,
such that
$\calc_t\not =\emptyset$ for each $t\in \R$. Further,  $S_I$ is bounded in
$W^{2,p}(\Omega)$, for each bounded  $I\subset \R$.
\item[2.]
Setting $\alpha_t=\inf\{\sup_\Omega u\,|\, u\in S_t\}$ and
$\overline{u}_t(x)=\sup\{u(x)\,|\, u\in S_t\}$, we have
$$
\lim_{t\to +\infty} \alpha_t=+\infty, \quad \mbox{and}\quad
\lim_{t\to -\infty}\sup_K \overline{u}_t(x)=-\infty,$$ for each
$K\subset\subset\Omega$, and  $\overline{u}_t<0$ in $\Omega$ for all $t$ below some number $t^-(h)$.
\end{itemize}
\end{teo}
The mere existence of
 solutions to \equ{equ1} when  $\lambda_1^-(F)\in (-L,0)$
was recently proved  in \cite{A}. Here we describe  qualitative properties
of the set of solutions.

To summarize,  the five theorems  above give a global picture
of the solutions of \re{equ1}, depending on the values of the eigenvalues
with respect to zero. This is shown on the following figure.
\begin{figure}[h]\label{figi}
\includegraphics[width=13.5cm, height=6cm]{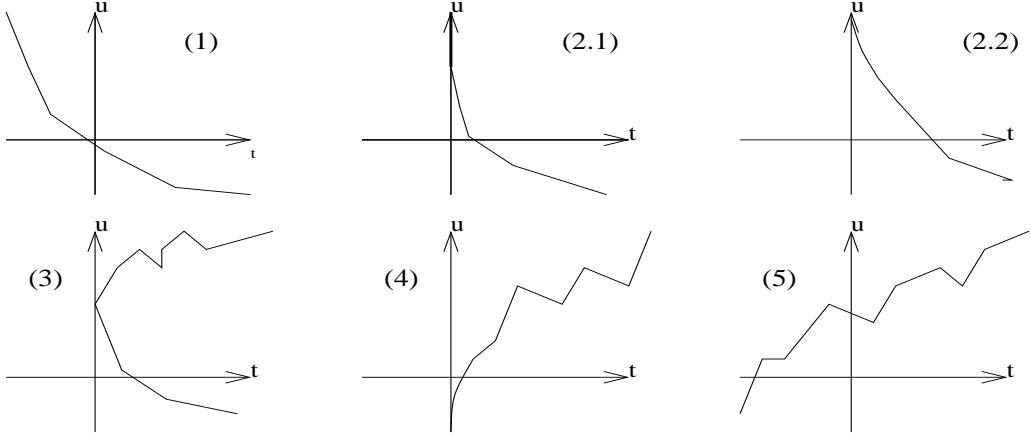}
\caption{\small\it The number at each graph corresponds to the number of the theorem where the shown situation is described.
 When $\lambda_1^+$ crosses $0$ the set $S$ curves so that
one region of non-existence and one region of multiplicity of solutions
 appears for $t$. Similarly when $\lambda_1^-$ crosses $0$
the set $S$
"uncurves" back.
 In this process, the set $S$ evolves
from being "decreasing", when both eigenvalues are positive,
 to being "increasing", at least for large $|t|$, when both eigenvalues are negative. Note (1) and (2.1)-(2.2) are exact, while in (3)-(5) there may be other solutions, except if Theorem \ref{teo6} below holds.}
\end{figure}

A number of remarks on  questions that are still open are now in order. First, it is clearly very important to give some characterization of the critical numbers $t^*$ in terms of $F, h,$ and $\lambda$.
On submitting this paper we learned of a very recent work by Armstrong \cite{A1}, where he studies this question in the case $\lambda=\lambda_1^+$, and proves part 1. of Theorem \ref{teo2}. by a different method. More specifically, he proves an interesting minimax formula for $\lpl(F)$, which generalizes the Donsker-Varadhan formula for linear operators to the nonlinear case. In particular, it is proved in \cite{A1} that
$$
\lpl = \min_{\mu\in \mathcal{M}(\Omega)}\sup_{u\in
C^2_+(\overline{\Omega})}\int_\Omega (-\frac{F(D^2u(x),Du(x),u(x),x)}{u(x)}\,d\mu(x).
$$
Further, if $\mathcal{M}^*$ is the subset of the set of probability measures $\mathcal{M}$, on which this minimum is attained, then for each $\mu\in \mathcal{M}^*$ there exists a positive function $\varphi_\mu\in L^{N/(N-1)}(\Omega)$ such that $d\mu = \varphi_\mu\varphi_1^+dx$, and  the number $t^*_+$ from Theorem 1.2 can be written as
$$
t^*_+ = -\min_{\mu \in\mathcal{M}^*} \int_\Omega h \varphi_\mu\,dx.
$$
 The  results in \cite{A1} and our Theorem 1.2 are complementary to each other, as we describe the set of solutions, while the main theorems in \cite{A1} characterize the critical value $t^*_+(h)$.

Next,
it is not clear how to distinguish between the two alternatives in
statement 4. of Theorem \ref{teo2} (that is, (2.1) and (2.2) on the above figure), for any given operator $F$. A simple and important example
 where we have (ii) is the Fucik operator
$F(u) = \Delta u + \lambda_1(\Delta)u^+ + bu^-$,  indeed if we had (i),
the fact that the solutions become positive for $t$ close to $t^*$, eliminates
the term in $u^-$, giving a contradiction. A rather simple example of an operator for which both (i) and (ii) can happen (depending on $f$) is given in \cite{A1}.

Naturally, the description of the  set $S$ when $\lambda_1^-=0$,
in contrast with  $\lambda_1^+=0$, is less precise due to the fact that
in this situation we only have degree theory  at our disposal to get existence
of solutions, and that uniqueness of solutions above $\lambda_1^-$ is not available in general (see however Theorem \ref{teo6} below).

Further, a number of basic questions can be asked about exact multiplicity
of solutions of \re{equ1} when $\lambda_1^+(F)<0$. When $\lambda_1^-(F)>0$
this question is a generalization of the famous Lazer-McKenna problem,
which concerns the Fucik equation
\beq\label{fucik}
 \Delta u+bu^+=\varphi_1\quad\mbox{ in }\; \Omega,\qquad
u=0\quad\mbox{ on }\; \partial\Omega.
\eqq
Here
 $F(D^2u,Du,u)= \Delta u + b u^+$,  $\lambda_1^+(F) =\lambda_1-b,
\lambda_1^-(F)=\lambda_1$, $b=\lambda_1^--\lambda_1^+$ and
 $\lambda_i$ are the eigenvalues of the Laplacian.
It is known that equation \re{fucik}
has exactly one solution if $b<\lambda_1$, exactly two solutions if
$b \in (\lambda_1,\lambda_2)$, exactly four solutions if
$b \in (\lambda_2,\lambda_3)$ and exactly six solutions if $b \in (\lambda_3,
\lambda_3+\delta)$,
see \cite{So} and the references in that paper.
This example suggests that multiplicity of solutions when the two
eigenvalues have opposite signs depends on the distance
$\lambda_1^--\lambda_1^+$. We conjecture that there exists a number
$C_0$ such that if $\lambda_1^+(F)<0<\lambda_1^-(F)\le \lambda_1^+(F) +
C_0$ then problem \re{equ1} has exactly two solutions, one solution or no solution, depending on $f$.

In the same way it should be asked if uniqueness of solutions
holds when $\lambda_1^-(F)\in (-L,0)$, for some $L>0$.
In view of the discussion  above one might expect that the answer is
affirmative if the two eigenvalues are sufficiently close to each other.
This fact  constitutes our last main theorem.
\begin{teo}\label{teo6}
There exists a number $d_0>0$ such that if $$ -d_0\le\lambda_1^+(F) \le \lambda_1^-(F)<0, $$ then problem \re{equ1} has at most one solution.
\end{teo}
A consequence of this result is that if  both Theorems \ref{teo5} and \ref{teo6} hold, then the sets $\calc$ of solutions obtained in Theorems \ref{teo4} and \ref{teo5} are continuous curves, like in Theorems \ref{teo1} and \ref{teo2}. We remark that $d_0$ is the difference between $\lpl(F,\Omega^\prime)$ and $\lpl(F,\Omega)$, where $\Omega^\prime$ is some subset of $\Omega$, whose Lebesgue measure is smaller than half the measure of $\Omega$ - see
Proposition \ref{teo8} and the proof of Theorem 1.6 in Section \ref{sect6}.

The article is organized as follows. In Section 2. we recall some
known results which we use repeatedly in our analysis. We also
complete the proof of Theorem \ref{teo1}.   Section \ref{sect3} is
devoted to   resonance phenomena at $\lambda_1^+=0$. In Section
\ref{sect4}  we analyze the existence and the properties of the set
of solutions of \equ{equ1} when $\lambda_1^-<0$. This set serves to
obtain the set of solutions at resonance when $\lambda_1^-=0$, in
Section \ref{sect5}. Finally, in Section \ref{sect6} we prove
Theorem~\ref{teo6}.

Some notational conventions will be helpful in the sequel.  When no
confusion arises, we write $F[u]:= F(D^2u, Du, u, x)$. We reserve
the notation
 $\norm{\cdot}=\norm{\cdot}_{\linf}$, while for all other norms we make
precise
mention to the corresponding space.

\setcounter{equation}{0}

\section{Preliminaries}\label{sect2}

In this section we give, for the reader's convenience, some of the results of the general theory of viscosity solutions of HJB equations, which we use in the sequel.
We start by restating the basic assumptions on the operator
$F:S_N\times \R^N\times\R\times\Omega\to\R$.
\begin{itemize}
\item[(H0)]  $F$ is
positively homogeneous of degree $1$, that is, for all $t\ge 0$ and
for all $(M,p,u,x)\in S_N\times \R^N\times\R\times\Omega$,
$$
F(tM,tp,tu,x)=tF(M,p,u,x).
$$
\item[(H1)] There exist $\gamma,\delta>0$ such that for all
$M,N\in S_N$, $p,q\in\R^N$, $u,v\in\R$, and a.e. $x\in\Omega$
\begin{eqnarray}
\M_{\lambda,\Lambda}^-(M-N)-\gamma|p-q|-\delta|u-v|
\le F(M,p,u,x)
 -F(N,q,v,x)\nonumber\\ ~~~~ \le
\M_{\lambda,\Lambda}^+(M-N)+\gamma|p-q|+\delta|u-v|.\nonumber
\end{eqnarray}
\item[(H2)]
$F(M,0,0,x)$ is continuous in $S_N\times\overline{\Omega}$.
\item[(H3)] If we denote $G(M,p,u,x)=-F(-M,-p,-u,x)$ then
\begin{eqnarray}
G(M-N,p-q,u-v,x)&\le &F(M,p,u,x)-F(N,q,v,x)\nonumber\\
&\le & F(M-N,p-q,u-v,x).\nonumber
\end{eqnarray}
\end{itemize}
Under (H0) the last assumption (H3) is equivalent to the convexity
of $F$ in $(M,p,u)$. The simple proof of this fact can be found for
instance in Lemma 1.1 in~\cite{QS2}. We recall that Pucci's extremal
operators are defined by $ {\cal M}^+(M)= \mathop{\sup}_{A\in {\cal
A}}\mbox{tr}(AM)$, $\mm(M)=  \mathop{\inf}_{A\in {\cal
A}}\mbox{tr}(AM), $ where $ {\cal A} \subset {\cal S}_N$ denotes the
set of matrices whose eigenvalues lie in the interval $[ \lambda,
\Lambda ]$.

We often use the following results from \cite{QS2} (Theorems 1.2-1.4 of that paper), which state that the principle eigenvalues are simple and isolated.

\begin{teo}[\cite{QS2}]\label{simple} Assume $F$ satisfies $(H0)-(H3)$ and there exists a viscosity solution
$u\in C(\overline{\Omega})$ of
\begin{equation}\label{simp1}
   F(D^2u, Du, u, x) =  -\lpl u\quad\mbox{in }\;\Omega,\qquad
u=0\quad\mbox{on }\;\partial\Omega,
\end{equation}
or of one of the problems
\begin{equation}\label{simp3}
\left\{ \begin{array}{rclll} F(D^2u, Du, u, x) &\le & -\lpl u&\mbox{in}&\Omega\\ u&>&0&\mbox{in}&\Omega,\end{array}\right.
\end{equation}
\begin{equation}\label{simp2}
\left\{ \begin{array}{rclll} F(D^2u, Du, u, x) &\ge & -\lpl u&\mbox{in}&\Omega\\
u(x_0)>0,\qquad u&\le&0&\mbox{on}&\partial\Omega,\end{array}\right.\end{equation}
for some $x_0\in \Omega$. Then $u\equiv t\fipl$, for some $t\in
\mathbb{R}$. If a function $v\in C(\overline{\Omega})$ satisfies
either (\ref{simp1}) or the inverse inequalities in (\ref{simp3}) or (\ref{simp2}),
with $\lpl$ replaced by $\lmin$, then $v\equiv t\fimin$ for some
$t\in \mathbb{R}$.
\end{teo}

\begin{teo}[\cite{QS2}]\label{isolated} There exists $\varepsilon_0>0$
depending on $N, \lambda, \Lambda, \gamma, \delta, \Omega,$ such
that the problem
\begin{equation}\label{sim}
   F(D^2u, Du, u, x) =  -\lambda u\quad\mbox{in }\;\Omega,\qquad
u=0\quad\mbox{on }\;\partial\Omega,
\end{equation}
has no solutions $u\not\equiv0$, for $\lambda\in(-\infty,
\lmin+\varepsilon_0)\setminus\{\lpl,\lmin\}$.
\end{teo}

 In the sequel we shall need the following {\it one-sided} Alexandrov-Bakelman-Pucci (ABP)
estimate, obtained in \cite{QS2} as well. The ABP inequality
for proper operators can be
found in \cite{CCKS} (an ABP inequality for the Pucci operator was
first proved in \cite{C}). We recall that $\lpl,\lmin$ are bounded above and below by constants which depend only on $N,\lambda,\Lambda,\gamma,\delta,\Omega,$ and that  both principal
eigenvalues of any proper operator are positive, see \cite{QS2}.
\begin{teo}[\cite{QS2}]\label{abp}
 Suppose the operator $F$ satisfies $(H0)-(H3)$.

{\bf I.} \ If $\lminf>0$ then for any $u\in C(\overline{\Omega})$,
$f\in L^N(\Omega)$, the inequality $F(D^2u, Du, u, x) \le f$
 implies
$$
\displaystyle \mathop{\sup}_{\Omega} u^- \le C(
\mathop{\sup}_{\partial\Omega} u^- + \norm{f^+}_{L^N(\Omega)}),\;
$$
where $C$ depends on  $N,\lambda, \Lambda, \gamma,
 \delta$, $\Omega$, and $1/\lmin$.

{\bf II.} In addition, if  $\lplf>0$ then  $F(D^2u, Du, u, x) \ge
f$ implies $$ \displaystyle \mathop{\sup}_{\Omega} u \le C(
\mathop{\sup}_{\partial\Omega} u^+ +\norm{f^-}_{L^N(\Omega)}). $$
Hence if  $\lplf>0$ then the comparison principle holds : if $u,v\in C(\overline{\Omega})$ are such that $F[u]\le F[v]$ in $\Omega$, $u\ge v$ on $\partial \Omega$, and one of $u,v$ is in $W^{2,N}(\Omega)$ then $u\ge v$ in $\Omega$.
\end{teo}

Note this result with $f=0$ gives {\it one-sided maximum principles}.  We also recall the following strong maximum principle  or Hopf's lemma, which is a
consequence from the results in \cite{BD} (a simple proof can be found in the appendix of \cite{A}).
\begin{teo}[\cite{BD}]\label{hopf}  Suppose $w\in C(\overline{\Omega})$
   is a viscosity solution of
$$
      \mm(D^2w) -\gamma|Dw|- \delta w\le 0\mbox{ in } \Omega,$$
      and $w\ge 0$ in $\Omega$.
 Then either $w\equiv 0$ in $\Omega$ or $w>0$ in $\Omega$ and
at any point $x_0\in\partial \Omega$ at which $w(x_0)=0$ we have $
\displaystyle\mathop{\liminf}_{t\searrow 0} \frac{w(x_0+t\nu)
-w(x_0)}{t} > 0, $ where $\nu$ is  the interior normal to
$\partial \Omega$ at $x_0$.
\end{teo}

We are going to use the following regularity result. It was proved in this generality in \cite{Sw} (interior estimate) and in \cite{W} (global estimate).

\begin{teo}[\cite{Sw},\cite{W}]\label{regul}  Suppose the operator $F$ satisfies $(H0)-(H2)$ and $u$ is a viscosity solution of $F(D^2u, Du, u, x) = f$ in $\Omega$, $u=0$ on $\partial \Omega$. Then $u\in W^{2,p}(\Omega)$, and
$$
\|u\|_{W^{2,p}(\Omega)} \le C\left(\|u\|_{\linf} + \|f\|_{L^p(\Omega)}\right),$$
where $C$ depends only on $N, p,\lambda, \Lambda, \gamma, \delta$, $\Omega$. \end{teo}

 Next we quote  the existence result from \cite{L} and \cite{QS2}.

\begin{teo}[\cite{QS2}]\label{existpos}  Suppose the operator $F$ satisfies $(H0)-(H3)$.

{\bf I.} \ If $\lmin(F,\Omega)>0$ then for any  $f\in
L^p(\Omega),p\ge N$, such that $f\ge 0$ in $\Omega$, there exists
a  solution $u\in W^{2,p}(\Omega)$ of $F(D^2u, Du, u, x) = f$ in $\Omega$,
$u=0$ on $\partial \Omega$, such that $u\le0$ in $\Omega$.

{\bf II.} In addition, if  $\lplf>0$ then  for any  $f\in
L^p(\Omega),p\ge N$, there exists a unique viscosity solution $ u\in W^{2,p}(\Omega)$ of
$F(D^2u, Du, u, x) = f$ in $\Omega$, $u=0$ on $\partial \Omega$.
\end{teo}

The next theorem is a simple consequence of the compact embedding $W^{2,p}(\Omega)\hookrightarrow C^{1,\alpha}(\Omega)$, Theorem \ref{regul}, and the convergence properties of viscosity solutions (see Theorem 3.8 in \cite{CCKS}).

\begin{teo}\label{conv}  Let $\lambda_n\to \lambda$ in $\mathbb{R}$ and $f_n\to f$ in $L^p(\Omega)$, $p>N$. Suppose  $F$ satisfies $(H1)$ and $u_n$ is a solution of $F(D^2u_n, Du_n, u_n, x) +\lambda_n u_n= f_n$ in~$\Omega$, $u_n=0$ on $\partial \Omega$. If $\{u_n\}$ is bounded in $L^\infty(\Omega)$ then a subsequence of $\{u_n\}$ converges in  $C^{1}(\overline{\Omega})$ to a function $u$, which solves $F(D^2u, Du, u, x) +\lambda u = f$ in $\Omega$,
$u=0$ on $\partial \Omega$. \end{teo}

We  now give the proof of Theorem \ref{teo1}.

\noindent{\bf Proof of Theorem \ref{teo1}}. Part 1. is a  consequence of Theorems \ref{abp} and \ref{existpos}.

Let us prove Part 2. For $t\in \R$, let $u_t$ be the solution
of \equ{equ1} with $f$ as in \equ{f}, that is, $F[u_t] = t\fipl +h$, where $\lpl(F)>0$. Then $\norm{u_t}/t$ is bounded as $t\to -\infty$. Indeed, if this is not the case, there exists
a sequence $\{t_n\}$ such that we have $t_n\to -\infty$ and  $\|u_{t_n}/t_n\|
\to \infty$, in particular $\|u_{t_n}\|\to \infty$.
Defining $\hat u_n= u_{t_n}/\|u_{t_n}\|$, we get by (H0)
$$
F(D^2 \hat u_n, D \hat u_n,\hat u_n,x)  =
\frac{t_n}{\|u_{t_n}\|}\varphi_1^++ \frac{h}{\|u_{t_n}\|}
\quad\mbox{in }\; \Omega,\qquad
\hat u_n = 0 \quad\mbox{on }\; \partial \Omega.
$$
The right-hand side of this equation converges to zero in $L^p(\Omega)$, so $\hat{u}_n$ converges uniformly to zero, by Theorem \ref{conv} (note the limit equation $F[\hat{u}]=0$ has only the trivial solution, since $\lplf>0$). This contradicts $\|\hat{u}_n\|= 1$.

Thus, by Theorem \ref{conv}, for some sequence $t_n\to -\infty$, we have that
\begin{equation}\label{rrr}\lim_{n\to \infty}{\frac{u_{t_n}}{-t_n}}=v^*
\quad\mbox{ in }\quad
C^1(\overline{\Omega}),
\end{equation}
 where $v^*$ satisfies
$$
F(D^2v^*, Dv^*,v^*,x)=-\varphi_1^+\quad\mbox{in}\quad \Omega,\qquad v^*=0 \quad\mbox{on}\quad \partial \Omega.
$$
By Theorems \ref{abp} and \ref{hopf} we have  $v^*>0$ in $\Omega$ and $\frac{\partial v}{\partial \nu}>0$ on $\partial \Omega$. These facts, \re{rrr}, and the monotonicity of $u_t$ in $t$ imply the last two statements of Part 2 (the analysis for $t\to \infty$ is similar).

That $S$ is Lipschitz follows from (H3) and Theorem \ref{abp}, applied to
$$
F[u_t-u_s]\ge (t-s)\fipl \qquad \mbox{ and } \qquad F[u_s-u_t]\ge (s-t)\fipl.
$$

Finally, the convexity property of the curve is a consequence of the following simple lemma and the comparison principle, Theorem \ref{abp}.
\begin{lema}\label{convexx}
Let $t_0,t_1\in \mathbb{R}$ and $t_k = kt_1 + (1-k)t_0$, for $k\in[0,1]$. Suppose $u_{t_i}\in S_{t_i}$, $i=0,1$. Then the function $ku_{t_1} + (1-k)u_{t_0}$ is a supersolution of
$$
   F(D^2u, Du, u, x) =  t_k\fipl +h \quad\mbox{in }\;\Omega,\qquad
u=0\quad\mbox{on }\;\partial\Omega.
$$
\end{lema}
\noindent{\bf Proof.} Use $F[ku_{t_1} + (1-k)u_{t_0}]\le kF[u_{t_1}] + (1-k)F[u_{t_0}]$.\hfill $\Box$

\bigskip

\noindent {\bf \underline{Notation}.} In what follows it will be convenient for us to write problem \re{first} in the form
\beq\label{princ}
\left\{\begin{array}{rclrr} F(D^2u,Du,u,x)+\lambda u &=&t\varphi_1(x) + h(x)&\mbox{in}& \Omega,\\
u&=&0&\mbox{on}& \Omega, \end{array}\right.\eqq
where $F$ is supposed to  proper (if necessary, we replace $F$ by $F-\delta$ and $\lambda$ by $\lambda + \delta$), and study its solvability in terms of the value of the parameter $\lambda\in \mathbb{R}^+$. For instance, Theorem \ref{teo2} corresponds to $\lambda= \lpl$, Theorem \ref{teo4} corresponds to $\lambda =\lambda_1^-$, Theorem \ref{teo1} corresponds to $\lambda <\lambda_1^+$,  etc.

\section{Resonance at $\lambda=\lambda_1^+$. Proof of Theorem \ref{teo2}}\label{sect3}

We first set up some preliminaries. Let $\{\lambda_n\}$ be a
sequence such that $\lambda_n<\lambda_1^+$ for all $n$, and
$\lim_{n\to\infty}\lambda_n=\lambda_1^+$. We consider the problem
 $$
 F(D^2u,Du,u,x) + \lambda_n u = t\fipl + h \quad \mbox{in}\quad \Omega,\qquad u=0 \quad\mbox{on}
\quad \partial \Omega,
$$
and  its unique solution   $u(n,t)$. In the sequel we shall write
$ u_n (t) = u(n,t) $ and also sometimes $u_n$ or $u_t$ instead of $u(n,t)$, when one of the parameters is kept fixed.

 We
define $\Gamma_n^+=\{u_n(
t)\,|\, t\in\R\}$. Recall that, by Theorem
\ref{teo1}, if $s<t$ then $u_n(t)<u_n(s)$.

We parameterize $\Gamma_n^+$ in the following way. We take a reference function $\tilde u_n=u_n(\tilde t_n)\in \Gamma_n^+$,  which is arbitrary but fixed for each $n\in \N$ (later we choose an appropriate sequence $\{ \tilde u_n\}$), and
we define the function
\begin{equation}\label{dn}
\left\{\begin{array}{rcl} d_n&:&\Gamma_n^+\to \mathbb{R}\\ d_n(u)&=&{\rm sign}{(u-\tilde u_n)}\|u-\tilde u_n\|.\end{array}\right.
\end{equation}
\begin{lema}\label{lema31}
The function ${d}_n :\Gamma_n^+\to\R$  is a bijection, for each $n\in\mathbb{N}$. In addition, $d_n$ is (Lipschitz) continuous.
\end{lema}
\noindent
{\bf Proof.}  By $(H3)$ for any $t_1,t_2\in \mathbb{R}$ (say $t_1>t_2$) we have
\beeq\label{diff1}
F[u_n(t_1) - u_n(t_2)]+ \lambda_n (u_n(t_1) - u_n(t_2)) \ge (t_1-t_2)\fipl.
\eeq
The ABP inequality (Theorem \ref{abp}) applies to this inequality - here we use  $\lambda_n<\lpl$ - so we have
$$
\|u_n(t_1) - u_n(t_2)\|\le C_n |t_1-t_2|.
$$
 If $t_1>t_2>\tilde t_n$ (the argument is the same if $t_2<t_1<\tilde t_n$) we get $$
|d_n(u_1) -d_n(u_2)| = \|u_{t_1} - \tilde u_n\| - \|u_{t_2} - \tilde u_n\| \le \|u_{t_1} - u_{t_2}\|\le C_n |t_1-t_2|.
$$
If $t_1>\tilde t_n>t_2$ we have
$$
|d_n(u_1) -d_n(u_2)| \le  \|u_{t_1} - \tilde u_n\| + \|u_{t_2} - \tilde u_n\| \le C_n (t_1-\tilde t_n+\tilde t_n-t_2) = C_n |t_1-t_2|,
$$
which proves the Lipschitz continuity.

Assume that $ d_n(u_n(t_1))= d_n(u_n(t_2))$, then $\|u_n(t_1)-\tilde u_n\|=
\|u_n(t_2)-\tilde u_n\|$ and $u_n(t_i)>\tilde u_n$ (or $u_n(t_i)<\tilde u_n$)
for $i=1,2$. On
the other hand, if $t_1\not=t_2$, say $t_1<t_2$, then $u_n(t_1)>u_n(t_2)$ and consequently
$\|u_n(t_1)-\tilde{u}_n\|\not=\|u_n(t_2)-\tilde{u}_n\|$, which
is impossible. Thus, $d_n$ is one to one.
By Part 2. in Theorem \ref{teo1} we see that ${d}_n$ is onto.
\hfill $\Box$

\medskip

Now we start the analysis of the resonance at $\lambda=\lambda_1^+$ (recall  we are working  with \re{princ}).
Given $s\in\R$ we define the proposition ${\cal P}(s)$ as follows:

\begin{itemize}
\item[${\cal P}(s):$]
{\it
There exist  sequences
$\{\lambda_n\}$, $\{h_n\}$ and $\{u_n\}$ such that $\lambda_n<\lambda_1^+$
for all~$n$,  $\lim_{n\to\infty}\lambda_n=\lambda_1^+$, $h_n\to h$ in $\lp$ as $n\to\infty$,
\begin{equation}\label{rety}
F(D^2u_n,Du_n,u_n,x)+\lambda_n u_n=s\varphi_1^++h_n,
\end{equation}
and  $\norm{u_n}$ is unbounded.}
\end{itemize}
By dividing \re{rety} by $\norm{u_n}$, thanks to $(H0)$, Theorem
\ref{simple} and Theorem \ref{conv}, we easily see that this
definition is equivalent to
\begin{itemize}
\item[${\cal P}(s):$]
{\it
There exist  sequences
$\{\lambda_n\}$, $\{h_n\}$, such that $\lambda_n<\lambda_1^+$
for all $n$, $\lambda_n\to\!\lambda_1^+$, $h_n\to h$ in $\lp$, the solution of
$
F(D^2u_n,Du_n,u_n,x)+\lambda_n u_n=s\varphi_1^++h_n
$
satisfies $\|u_n\|\to \infty$, and  $$\frac{u_n}{\|u_n\|}\to \fipl>0\qquad \mbox{in }\; C^1(\overline{\Omega}).$$}
\end{itemize}
We define
\begin{equation}
t^*_+=\sup\{t\in \R\,|\, {\cal P}(s) \mbox{ for all } s<t\}.
\end{equation}
The next lemmas give meaning to this definition.
\begin{lema} Given $\bar t\in \mathbb{R}$,
${\cal P}(\bar t)$ implies
${\cal P}( t)$ for all $t<\bar t$.
\end{lema}
\noindent
{\bf Proof.}
Assuming the contrary,
there is $t_0<\bar t$ such that ${\cal P}(t_0)$ is false.
This means that for some sequences $\{\lambda_n\}$, $\{h_n\}$ as above, the sequence of the solutions of
$$
F(D^2v_n, Dv_n,v_n,x)+\lambda_n v_n=\bar t\varphi_1^++h_n\quad \mbox{in}\quad \Omega,\qquad v_n=0 \quad\mbox{on}
\quad \partial \Omega.
$$
is unbounded, while the sequence of the  solutions of
$$
F(D^2u_n, Du_n,u_n,x)+\lambda_n u_n=t_0\varphi_1^++h_n\quad \mbox{in}\quad \Omega,\qquad u_n=0 \quad\mbox{on}
\quad \partial \Omega,
$$
is bounded in $\linf$. By the comparison principle (Theorem \ref{abp})
$v_n\le u_n$ for all $n$, since $\bar t>t_0$ and $\fipl>0$. On the other hand,
by the one-sided ABP inequality, Theorem \ref{abp} 1. (note $\lambda_n$ is uniformly away from $\lambda_1^-$, that is, $\lambda_1^-(F+\lambda_n)\ge \lmin(F)-\lpl(F)>0$), the sequence $\{v_n\}$ is bounded below.
Thus $\{v_n\}$ is bounded, a contradiction. \hfill $\Box$

\begin{lema}\label{lem33}
There exists a real number $\bar t_1 = \bar t_1(h)$, such that the problem
\beeq\label{rref}
F(D^2u, Du,u,x)+\lambda_1^+ u=t\varphi_1^++h\quad \mbox{in}\quad \Omega,
\qquad u=0 \quad\mbox{on}
\quad \partial \Omega,
\eeq
has no solutions for $t< \bar t_1$.
\end{lema}
\noindent
{\bf Proof.} Let $v$ be the solution of the Dirichlet problem
$$
F(D^2v, Dv,v,x)=-h\quad \mbox{in}\quad \Omega,
\qquad u=0 \quad\mbox{on}
\quad \partial \Omega
$$
(this problem is uniquely solvable, by the well-known results on proper equations, or by Theorem \ref{existpos}). We are going to show that the statement of the lemma is true with $$\bar t_1 = -1- \lpl \sup_{x\in \Omega}\frac{v(x)}{\fipl(x)}.$$
The last quantity is finite, by Theorems \ref{abp}-\ref{regul}.

Indeed, if  \re{rref} has a solution $u=u(t)$ for some  $t<\bar t_1$, we get
\begin{eqnarray}
F[u+v] + \lpl(u + v) &\le& F[u] + F[v] +\lpl u + \lpl v\nonumber\\
&\le & t\fipl + \lpl v\le -\fipl<0\label{reftr},
\end{eqnarray}
where we have used $F[u+v]\le F[u]+ F[v]$, which follows from (H3). Since we have  $\lmin(F+\lpl, \Omega)= \lmin - \lpl>0$, Theorem \ref{abp} 1. again applies and yields $u+v>0$ in $\Omega$. We can now use Theorem \ref{simple}  and conclude that $u + v$ is a multiple of $\fipl$, which contradicts the strict inequality in \re{reftr}. \hfill $\Box$

\begin{lema}\label{lemr}
The set $T=\{t\in \R\,|\, {\cal P}(t)\}$ is not empty.
\end{lema}
\noindent
{\bf Proof.}
Assuming the contrary, we find sequences $\{t_n\}$, $\{u_n^m\}$, such
that ${\cal P}(t_n)$ is false, $t_n\to-\infty$ as $n\to\infty$,  $u_n^m$ satisfies
$$
F(D^2u_n^m, Du_n^m,u_n^m,x)+(\lambda_1^+-1/m) u_n^m=t_n\varphi_1^++h\quad \mbox{in}\quad \Omega,
\,\,  u_n^m=0 \quad\mbox{on}
\quad \partial \Omega,
$$
for each $n$, and $\{u_n^m\}$ is bounded in $\linf$ as $m\to \infty$. Hence, by Theorem~\ref{conv}, $u_n^m$ converges as $m\to \infty$ (up to a subsequence), for each fixed $n$,  to a function $u_n$ which satisfies \re{rref} with $t=t_n$. This and the previous lemma give a contradiction, when $t_n$ is sufficiently small. \hfill $\Box$

\begin{lema}\label{lem35}
The set $T$ is bounded above, that is, $t^*_+$ is a real number.
\end{lema}
\noindent
{\bf Proof.}
Let $\lambda_n\nearrow\lambda_1^+$, $h_n\to h$ in $\lp$, and let $u_n=u_n(t)$ be such that
$$
F(D^2u_n,Du_n,u_n,x)+\lambda_n u_n=t\varphi_1^++h_n\quad \mbox{in}\quad \Omega,
\qquad u_n=0 \quad\mbox{on} \quad \partial \Omega.
$$
(we recall this problem has a unique solution, since $\lambda_n<\lambda_1^+$ and comparison holds). We need to show $\{u_n\}$ is bounded in $\linf$, if $t$ is large enough.

First,  Theorem \ref{abp} I. implies that $u_n$ is bounded below independently of~$n$ (we recall once again that $\lmin(F+\lambda_n)\ge \lmin-\lpl>0$).

Next, let
$v_n$ be the solution of the Dirichlet problem
$$
F(D^2v_n,Dv_n,v_n,x)=\min\{h_n,0\}\le 0\quad\mbox{in}\quad\Omega, \qquad v=0
\quad\mbox{on} \quad \partial \Omega.
$$
  Then $v_n\ge 0$   in $\Omega$, by the maximum principle, $\{v_n\}$ is bounded in $C^1(\overline\Omega)$, by Theorems \ref{abp} and \ref{regul}, and
$$
F[v_n]+\lambda_n v_n\le \min\{h_n,0\}+\lambda_1^+ v_n\le h_n+t\varphi_1^+= F[u_n]+\lambda_n u_n,
$$
provided
\begin{equation}\label{fkij} t>\lpl \sup_{x\in \Omega, n\in \mathbb{N}}\frac{v_n(x)}{\fipl(x)}\,.\end{equation}

By the comparison principle $u_n\le v_n$, hence $u_n$ is bounded above independently of $n$. So ${\cal P}(t)$ is false if \re{fkij} holds. \hfill $\Box$

\medskip

The following two propositions give existence and uniqueness of solutions
to our  problem at resonance, provided  $t>t^*_+$.
\begin{proposition}\label{lema33}
1. If $t>t^*_+$ then the equation
\begin{equation}\label{retz1}
F(D^2u, Du,u,x)+\lambda_1^+ u=t\varphi_1^++h\quad \mbox{in}\quad \Omega,\qquad u=0 \quad\mbox{on}
\quad \partial \Omega,
\end{equation}
possesses at least one solution.

2.  If $t<t^*_+$ then \re{retz1} has no solutions.
\end{proposition}
\noindent
{\bf Proof.} 1. Given a  sequence $\{\lambda_n\}$ such that
$\lambda_n<\lambda_1^+$ for all
$n\in\N$ and $\lambda_n\to\lambda_1^+$ as $n\to\infty$,
there is a sequence $\{u_n\}$ such that
\begin{equation}\label{retz2}
F(D^2u_n,Du_n,u_n,x)+\lambda_n u_n=t\varphi_1^++h\quad \mbox{in}\quad \Omega,\qquad u_n=0 \quad\mbox{on}
\quad \partial \Omega.
\end{equation}
Then $t>t^*$ implies $\{u_n\}$ is bounded, so by Theorem \ref{conv} $\{u_n\}$ converges, up to a subsequence, to a function $u$ satisfying \re{retz1}.

2. Suppose for contradiction \re{retz1} has a solution $u$ for some $t<t^*_+$. Fix $t_1\in (t,t^*_+)$. Then ${\cal P}(t_1)$ holds, so for some sequences $\lambda_n \nearrow \lambda_1^+$, $h_n\to h$,  the sequence of solutions $u_n$ of
$$
 F(D^2u_n,Du_n,u_n,x)+\lambda_n u_n=t_1\varphi_1^++h_n\quad \mbox{in}\quad \Omega,\qquad u_n=0 \quad\mbox{on}
\quad \partial \Omega
$$is such that $u_n\ge k_n\fipl$ for some $k_n\to \infty$. Let now $w_n$ be the solution of
\begin{equation}\label{retz3}
 F(D^2w_n,Dw_n,w_n,x)=h_n-h\quad \mbox{in}\quad \Omega,\qquad u_n=0 \quad\mbox{on}
\quad \partial \Omega
\end{equation}
By Theorems \ref{abp} and \ref{regul} we know that (up to a subsequence) $w_n\to 0$ in $C^1(\overline{\Omega})$.
Hence, by the boundary Lipschtiz estimates (see Theorem \ref{regul}, or Proposition 4.9 in \cite{QS2}) applied to \re{retz1}, \re{retz3},  we have
$$
\|u\| + \|w_n\|\le C\mbox{dist}(x,\partial\Omega),
$$
which implies
$$
u_n-w_n-u>0$$
 for $n$ sufficiently large. Since $w_n\to 0$ and $\lambda_n\to\lpl$ we also have
 $$ t_1\fipl -2\lambda_1^+|w_n|>t\fipl  \qquad \mbox{and} \qquad |u|\le \frac{t_1-t}{2(\lambda_1^+-\lambda_n)}\fipl\quad\mbox{in }\;\Omega. $$
However $(H3)$ implies $F[u_n-w_n-u]\ge F[u_n]-F[w_n]-F[u]$, so
$$
F[u_n-w_n-u] + \lambda_n(u_n-w_n-u)\ge (t_1-t)\fipl - \lambda_1^+|w_n|+ (\lambda_1^+-\lambda_n) u \ge 0.
$$
Then the maximum principle (Theorem \ref{abp}) gives $u_n-w_n-u\le 0$, a contradiction. \hfill $\Box$

\medskip

Next we prove the uniqueness of solutions above $t^*_+$. In order to do this,
we need the following simple result on convex functions.

\begin{lema}\label{convex}
Let $f:\R^n\to\R$ be positively homogeneous of degree one and convex.
If for some $u,v\in \R^n$ and for some $\tau>0$ we have
\begin{equation}\label{hyp}
f(u+\tau v)=f(u)+\tau f(v)
\end{equation}
then \re{hyp} holds for all $\tau\ge 0$.
\end{lema}
\noindent
{\bf Proof.}
Using \equ{hyp} and the homogeneity of $f$ we find that
$$
f(\lambda_0u+(1-\lambda_0)v)=\lambda_0 f(u)+(1-\lambda_0)f(v),
$$
with $\lambda_0=1/(1+\tau)$.
If there is $\lambda\in (\lambda_0,1)$ such that
\begin{equation}\label{i2}
f(\lambda u+(1-\lambda )v) < \lambda  f(u)+(1-\lambda )f(v)
\end{equation}
we take $\theta= 1-\lambda_0/\lambda\in (0,1)$, that is, $(1-\theta)\lambda=\lambda_0$, and notice that
\begin{eqnarray}
\lambda_0f(u)+(1-\lambda_0)f(v)&=&f(\lambda_0 u+(1-\lambda_0)v)\nonumber \\
&=& f\left( \theta v + (1-\theta)(\lambda u + (1-\lambda)v)\right)\nonumber\\
&\le &\theta f(v)+(1-\theta)f(\lambda u+(1-\lambda)v) \nonumber \\
&<& (1-\lambda_0)f(v)+\lambda_0f(u),\nonumber
\end{eqnarray}
which is a contradiction. If there is $\lambda \in (0,\lambda_0)$ such that \equ{i2} holds,
we proceed similarly. Thus,
$
f(\lambda u+(1-\lambda )v) = \lambda  f(u)+(1-\lambda )f(v),
$
for all $\lambda\in [0,1]$. From here we get the conclusion, taking
$\lambda=1/(1+t)$.\hfill  $\Box$

\begin{proposition}\label{prop1}
1. If $t>t^*_+$ and $u_1, u_2$ satisfy
$$
F(D^2 u_i,D u_i, u_i,x)+\lambda_1^+ u_i=t\varphi_1^++h
\quad \mbox{in}\quad \Omega,\quad u_i=0 \quad\mbox{on} \quad \partial \Omega,
$$
$i=1,2$, then $u_1=u_2$.

2.  If $t=t^*_+$ and $u_1, u_2$ are as in 1. then $u_1= u_2 + s\varphi_1$, for some $s\in\mathbb{R}$.
\end{proposition}
\noindent
{\bf Proof.}
Suppose $u_1\not =u_2$, then we may assume  there exists $x_0\in\Omega$
such that $u_1(x_0)>u_2(x_0)$. By (H3) we have
$
F[u_1-u_2]+\lambda_1^+ (u_1-u_2)\ge 0,
$
so by Theorem \ref{simple} there exists $\tau>0$ such that
$u_1-u_2=\tau\varphi_1^+$. This implies
\begin{equation}\label{equality}
F[u_1+\tau\varphi_1^+]=F[u_1]+\tau F[\varphi_1^+]\quad\mbox{a.e. }\;\mbox{in }\;\Omega
\end{equation}
(note $u_1,\fipl\in W^{2,N}(\Omega)$).
Consider the function $f(X)=F(X,x)$ where $X=(M,p,u)\in S_N\times\R^N\times\R=\R^{N^2+N+1}$, and $x\in\Omega$ is  fixed. By hypotheses
(H0) and (H3) the function $f$ is positively homogeneous of degree one and convex in $X$.
Therefore  we can use Lemma \ref{convex} to conclude that \equ{equality} holds for every $\tau>0$.

We obtain that for every $n\in\mathbb{N}$ the function
$v_n=u_1+n\varphi_1^+$ satisfies
$
F(D^2 v_n,D v_n, v_n,x)+\lambda_1^+ v_n=t\varphi_1^++h$ in $\Omega$, $u_i=0 $ on $\partial \Omega$. It follows that
$$
F[v_n] + \left(\lambda_1^+-\frac{1}{n^2}\right)v_n = t\varphi_1^++h - \frac{1}{n^2}u_1- \frac{1}{n}\fipl =: t\varphi_1^++h_n
$$
in $\Omega$. Note $h_n\to h$ in $L^p(\Omega)$. However this is impossible if $t>t^*_+$, by the definition of $t^*_+$, since $\norm{v_n}$ is unbounded, which means ${\cal P}(t)$ holds.
\hfill $\Box$

\medskip

Now we study the behaviour of the branch $\Gamma_n^+$ as $n\to\infty$. Let $\tilde u$
be the unique solution (given by Proposition \ref{lema33}) of
$$
F(D^2 \tilde u,D\tilde u,\tilde u,x)+\lambda_1^+ \tilde u=(1+t^*_+)\varphi_1^++h
\quad \mbox{in}\quad \Omega,\quad \tilde u=0 \quad\mbox{on} \quad \partial \Omega,
$$
and set
$$
d(u) ={\rm sign}{(u-\tilde u)}\|u-\tilde u\|.
$$
\begin{lema}\label{lema36}
If $u_i$ and $t_i$, $i=1,2$, are such that $d(u_1)=d(u_2)$ and
$$
F(D^2 u_i,D u_i, u_i,x)+\lambda_1^+ u_i=t_i\varphi_1^++h
\quad \mbox{in}\quad \Omega,\quad u_i=0 \quad\mbox{on} \quad \partial \Omega,
$$
for $i=1,2$, then $t_1=t_2$ and $u_1=u_2$.
\end{lema}
\noindent
{\bf Proof.} By Proposition \ref{prop1} $u_1\not=u_2$ implies $t_1\not=t_2$. If $t_1\not =t_2$ (say $t_1>t_2$),
$$
F[u_1-u_2]+\lambda_1^+ (u_1-u_2) \ge
 (t_1-t_2)\varphi_1^+>0
\quad \mbox{in}\quad \Omega,\quad u_1-u_2=0 \quad\mbox{on} \quad \partial \Omega.
$$
Either there exists $x_0\in\Omega$ such that $u_1(x_0)>u_2(x_0)$ or $u_1\le u_2$ in
$\Omega$. In the first case, Theorem \ref{simple}
implies the existence of $\tau>0$ such that
$u_1-u_2=\tau \varphi_1^+$ so that $u_1>u_2$ in $\Omega$.
In the second case, by the strong maximum principle, we have that $u_1=u_2$ (excluded by $t_1\not= t_2$) or
$u_1<u_2$ in $\Omega$.

Thus, if $u_1\not =u_2$, then either $u_1>u_2$ or $u_1<u_2$ in $\Omega$, and in both cases
$d(u_1)\not = d(u_2)$, completing the proof of the lemma.\hfill
$\Box$

\medskip

We recall (Lemma \ref{lema31}) that the set $\Gamma_n^+$
can be re-parameterized as
a curve, by using the function $d_n$. In the definition of $d_n$ we used the arbitrary function $\tilde u_n$, which we choose now as  the unique
solution of
$$
F(D^2 \tilde u_n,D\tilde u_n,\tilde u_n,x)+\lambda_n \tilde u_n=(1+t^*_+)\varphi_1^++h
\quad \mbox{in}\quad \Omega,\quad \tilde u_n=0 \quad\mbox{on} \quad \partial \Omega.
$$
By the definition of $t^*_+$ $\{\norm{\tilde u_n}\}$ is bounded, so by Theorem \ref{conv} and the uniqueness property
proved in Proposition
\ref{prop1} we find that $\tilde u_n\to \tilde u$, where $\tilde u$ is as above,  the unique solution of
$$
F(D^2 \tilde u,D\tilde u,\tilde u,x)+\lpl \tilde u=(1+t^*_+)\varphi_1^++h
\quad \mbox{in}\quad \Omega,\quad \tilde u=0 \quad\mbox{on} \quad \partial \Omega.
$$
By Lemma \ref{lema36}, for fixed
$d\in\R$ the following system in $(u,t)$
\begin{equation}
\left\{\begin{array}{rcl}F(D^2 u,D u, u,x)+\lambda_n u&=&t\varphi_1^++h
\quad \mbox{in}\quad \Omega,\quad u=0 \quad\mbox{on} \quad \partial \Omega,\\ d_n(u)&=&d
\end{array}\right.\end{equation}
has a unique solution $(u_n, t_n)$ in $C(\overline{\Omega})\times \mathbb{R}$. The sequence $\{u_n\}$ is bounded, since $\|u_n-\tilde{u}_n\|=|d|$ and $\{\tilde{u}_n\}$ is bounded. Hence  $\{t_n\}$ is also bounded (if not, $u_n/t_n \to 0$, so by passage to the limit $F[0]=\varphi_1^+$).

 Then  subsequences of $\{u_n\}$ and  $\{t_n\}$  converge to a function $u=u(d)$ and a number  $t= t(d)$, which satisfy
\begin{equation}\label{eq11}
F(D^2 u,D u, u,x)+\lambda_1^+ u=t\varphi_1^++h
\quad \mbox{in}\quad \Omega,\quad u=0 \quad\mbox{on} \quad \partial \Omega.
\end{equation} By Lemma
\ref{lema36} the whole
sequences $\{u_n\}$ and $\{t_n\}$ converge to the same limit that we
call $u(d)$  and $t(d)$.

\begin{lema} The map $\left\{\begin{array}{rcl} U : \mathbb{R} &\to& C(\overline{\Omega})\times \mathbb{R}\\ U(d) &=& \left( u(d),t(d)\right)\end{array}\right.$ is continuous.
\end{lema}
\noindent
{\bf Proof.} Take $d_k\to d$ as $k\to\infty$. Then the sequences $u_k=u(d_k)$, $t_k=t(d_k)$ are bounded, as above. Any convergent subsequence of $\{(u_k,t_k)\}$ tends to a solution of an equation to which Lemma \ref{lema36} applies, so the whole sequences $u_k$, $t_k$ converge to $u(d)$, $t(d)$. \hfill $\Box$
\medskip

We  define $\Gamma^+ = \{ u(d) \,|\, d\in\R\}$. The last lemma allows us to say that $\Gamma^+$ is actually a continuous curve, the pointwise limit of the
curves $\{\Gamma_n^+\}$.

\begin{lema}\label{sij} If $t_1>t_2\ge t_+^*$ then any two solutions of $$
F(D^2 u_i,D u_i, u_i,x)+\lambda_1^+ u_i=t_i\varphi_1^++h
\quad \mbox{in}\quad \Omega,\quad u_i=0 \quad\mbox{on} \quad \partial \Omega,
$$
are such that  $u_1<u_2$ in $\Omega$.
\end{lema}
\noindent
{\bf Proof.} We already showed in the proof of Lemma \ref{lema36} that either $u_1>u_2$ or $u_1<u_2$ in $\Omega$. Since the curve $\Gamma^+$ is the limit of $\Gamma_n^+$ which are strictly decreasing in $t$, $u_1>u_2$ is impossible. \hfill $\Box$

\medskip

\noindent
{\bf Proof of Theorem \ref{teo2}.} The set of solutions is
$\{(u(d),t(d))\,|\, d\in\R\}$, as the above discussion shows.
Part 1. of the Theorem was proved in Proposition~\ref{lema33}. The first two statements of Part 2. follow from Proposition \ref{prop1} and Lemma \ref{sij}.

For $t>t^*$, let $u_t$ be the solution of
\begin{equation}\label{tukat}
F(D^2 u_t,D u_t, u_t,x)+\lambda_1^+ u_t=t\varphi_1^++h
\quad \mbox{in}\quad \Omega,\quad u_t=0 \quad\mbox{on} \quad \partial \Omega,
\end{equation}
By Lemma \ref{sij} $u_t$ is strictly decreasing in $t$.

When $t\to t^*_+$ two cases may occur : either $\norm{u_t}$ is bounded or $\norm{u_t}\to \infty$. In the first case the monotonous sequence $u_t$ converges in  $C^1(\overline{\Omega})$ to a solution $u^*$ of \re{tukat} with $t= t^*_+$. Then by Proposition \ref{prop1} all solutions $u\in \Gamma^+$ with $d(u)\ge d(u^*)$ are solutions of  \re{tukat} with $t= t^*_+$, which is the situation described in Part 4. (ii). In the second case $u_t/\norm{u_t}$ converges in  $C^1(\overline{\Omega})$ to $\varphi_1^+>0$ which implies Part~4.~(i). Note in this case there cannot be solutions with $t=t^*_+$, because of Lemma~\ref{sij}.

Let us now consider the limit $t\to \infty$.
First, if for some sequence $t_n\to \infty$ we have $\norm{u_{t_n}}/{t_n}\to 0$, we divide \re{tukat} by $t_n$, pass to the limit and get a contradiction. So $\norm{u_t}\to \infty$ as $t\to \infty$. By the monotonicity of $u_t$ in $t$, $\min_\Omega u_t<-1$ for sufficiently large $t$.

Assume for some sequence $t_n\to \infty$ we have $\norm{u_{t_n}}/t_n\to \infty$. Then we divide \re{tukat} by $\norm{u_{t_n}}$ and see that $u_{t_n}/\norm{u_{t_n}} \rightrightarrows \varphi_1^+$, which is impossible, since $u_{t_n}$ takes negative values and $\varphi_1^+>0$.

So there is a sequence $t_n\to \infty$ such that $u_{t_n}/\norm{u_{t_n}}$ converges in $C^1(\overline{\Omega})$  to a solution of
\begin{equation}\label{getit}
F(D^2v, Dv, v,x) + \lambda_1^+ v = k\varphi_1^+\quad \mbox{ in }\; \Omega,\qquad v=0\quad\mbox{on }\,\partial \Omega,\end{equation}
 for some $k>0$. This  problem is the particular case of \re{princ} when $h=0$. It is clear that \re{getit} has solutions for $k\ge 0$ (by Theorem \ref{existpos}) and does not have solutions for $k<0$ (by the definition of $\lambda_1^+$ and Theorem \ref{simple}). Further, this problem obviously has solutions for $k=0$ (in other words, for $h=0$ we always are in the case 4. (ii)) and the minimal solution at $k=0$ is $u^*=0$. Then, by the properties  of the curve of solutions we already proved \re{getit} has a unique solution which satisfies $v<u^*=0$, since $k>0$.

This means $u_{t_n}/\norm{u_{t_n}}$ converges in $C^1(\overline{\Omega})$ to a negative function $v$, such that $\frac{\partial v}{\partial \nu}<0$ on $\partial \Omega$ (by \re{getit} and Hopf's lemma). This implies statement~3. for the subsequence $u_{t_n}$. Since $u_t$ is monotonous, we have 3. for all $u_t$.

Finally, let us show that $t\to u_t(x)$ is convex.  With the notations from Lemma \ref{convexx}, we note that for each compact interval $[t_0, t_1]\subset [t^*_+,\infty)$ there exists a function $v\in W^{2,p}(\Omega)$ which is a subsolution of
\begin{equation}\label{rezaq}
   F(D^2u, Du, u, x) +\lpl u =  t_k\fipl +h \quad\mbox{in }\;\Omega,\qquad
u=0\quad\mbox{on }\;\partial\Omega,
\end{equation}
and $v< ku_{t_1}+ (1-k)u_{t_0}$, for each $k\in(0,1)$ (we take $u_{t_0}=u^*$ if $t_0=t^*_+$). For instance, we can take $v$ to be the negative solution -- given by Theorem \ref{existpos} I. -- of the problem
$$
   F[v] +\lpl v =  \max\{t_1,1\}\fipl +\max\{h,0\} \quad\mbox{in }\;\Omega,\qquad
u=0\quad\mbox{on }\;\partial\Omega,
$$
and then a take a multiple of  $v$ by a sufficiently large constant, to ensure that  $v<u_{t_1}\le ku_{t_1}+ (1-k)u_{t_0}$ for each $k\in(0,1)$. Then by Lemma \ref{convexx} and the usual sub- and supersolution method there exists a solution of  \re{rezaq} which is below  $ku_{t_1}+ (1-k)u_{t_0}$. By the uniqueness which we already proved, this solution has to be $u_{t_k}$.

Theorem \ref{teo2} is proved. \hfill
$\Box$

\setcounter{equation}{0}
\section{The case $\lambda>\lambda_1^-$.
Proof of Theorem~\ref{teo5}}\label{sect4}

In this section we  prove
Theorem \ref{teo5} and some auxiliary results which will be helpful in our
analysis of the resonance phenomena at $\lambda=\lambda_1^-$.

We start with some simple preliminary lemmas which will lead us to the proof of the first part in Theorem \ref{teo5}. Our arguments for Lemmas \ref{lemche}-\ref{lemche1} below are
similar to those in \cite{BEQ},\cite{QS2} and  \cite{A}, but we  sketch them here
for completeness.
We define the operators
$$F_\s (D^2  u,D  u, u,x)=\s F(D^2  u,D  u, u,x)+(1-\s )\Delta u$$
and we write $\lambda_1^-(\s)=\lambda_1^-(F_\s)$, for $\s\in [0,1]$.
Note that $F_\tau$ satisfies $(H0)-(H3)$ and, recalling that we work with
 \re{princ}, $F_\tau$ is proper, since $F$ is
proper.
\begin{lemma}\label{lemche}   The function $\s\to\lambda_1^-(\s)$ is continuous  in
the interval $ [0,1]$ and  there exists $\bar{\varepsilon}>0$ so
that there is no eigenvalue of $F_\s$
in the interval  $(\lambda_1^-(\s),\lambda_1^-(\s)+\bar{\varepsilon}]$, for $\s\in [0,1]$.
\end{lemma}

\noindent
{\bf Proof.}  Let $\{\s_n\}$ be a sequence in $[0,1]$, then it follows by
Proposition 4.1 in \cite{QS2} that the sequence $\{
\lambda_1^-(\s_n)\}$  is bounded. Then, by a compactness argument and
the simplicity of the  eigenvalues, the continuity follows.
The isolation property follows by the same argument as the one used in the proof of
Theorem 1.3 in  \cite{QS2}.\hfill  $\Box$

\begin{lemma}\label{lemche1} There exists $\varepsilon>0$ such that          for each $\lambda \in(\lambda_1^-, \lambda_1^-+\varepsilon)$ and each  $n\in \mathbb{N}$ there is a closed
connected set
$\calc(\lambda,n) \subset C(\overline\Omega)\times [-n,n]$, with the property  that for all
$(u,t)\in
\calc(\lambda,n)$ we have
$$
F(D^2u,Du,u,x)+\lambda u= t \varphi_1^++h
\quad \mbox{in}\quad \Omega,\quad  u=0 \quad\mbox{on} \quad \partial \Omega.
$$
Moreover, if we define the projection
${\cal P}:  C(\overline\Omega)\times \RR\to \RR$ as
${\cal P}(u,t)=t$, we have
${\cal P}(\calc(\lambda,n))=[-n,n]$. \end{lemma}

\noindent
{\bf Proof.} For $\s\in [0,1]$, let us define
$$ \lambda_2(\s)=\inf\{\mu> \lambda_1^-(\s) \,|\, \mu\; \mbox{ is an
eigenvalue of }\; F_\s\,\mbox{ in } \Omega \}.$$
Observe that $\lambda_2(\tau)=+\infty$ is possible.
Then, given $\lambda\in (\lambda_1^-, \lambda_1^-+\bar \varepsilon)$,
by the previous lemma there exists a continuous function $\mu:[0,1] \to \RR$
such that
$\mu(1)=\lambda$, $\lambda_1^-(\tau)<\mu(\tau)<\lambda_2(\tau)$ and the equation
\begin{equation}\label{equei}
 F_\tau (D^2u,Du,u,x)+\mu(\tau) u=0\quad\mbox{ in }\; \Omega,\qquad
u=0\quad\mbox{ on }\; \partial\Omega, \end{equation}
has no non-trivial solution, for all $\tau\in [0,1]$.
Now we define the operator $G:\RR\times[0,1] \times C(\overline\Omega)\to
C(\overline\Omega)$, for $(t, \tau, v) \in \RR\times[0,1] \times C(\overline\Omega)$
as  $u=G(t,\tau,v)$, where $u$ is  the solution of the equation
\begin{equation}\label{pr1} F_\tau(D^2u,Du,u,x)=-\mu(\tau) v+t\varphi_1^++h \quad\mbox{ in }\; \Omega,\qquad
u=0\quad\mbox{ on }\; \partial\Omega. \end{equation}
When we restrict the variable $t$ to the interval $[-n,n]$, the
operator $G$  becomes
compact.
Moreover, there exists $R>0$ such that the Leray-Schauder degree
$d(I-G(t,\tau,.),B_R,0)$ is well defined. Indeed, a priori bounds follow
directly from the non-existence property of equation \equ{equei}, in fact,
 if \re{pr1} has a sequence of solutions $u_n=v_n$ with
$\norm{u_n}\to\infty$ we divide \re{pr1} by $\norm{u_n}$, pass to the limit
and get a contradiction.
Then, by the homotopy invariance of the Leray -Schauder degree, we have
$$d(I-G(t,1,\cdot),B_R,0)=d(I-G(t,0,\cdot),B_R,0)=-1.$$
The  last equality is a standard fact,  since the operator $F_0$ is the
Laplacian.
Thus, by the well-known results  of \cite{R} (alternatively, we refer to \cite{chang}), see in particular Corollary 10
in chapter V of that work, the lemma follows. \hfill $\Box$

\medskip

We will need the following topological result, whose proof is a direct consequence of
Lemma 3.5.2 in \cite{chang}.
\begin{lemma}\label{conn}
Let  $R\subset C(\overline\Omega)\times [-n,n]$ be a compact connected set
such that
${\cal P}(R)=[-n,n]$. If $R_0=\{(u,t)\in R\,
|\, t\in [t_-,t_+]\}$, with $[t_-,t_+]\subset
[-n,n]$
then there exists a connected component $E_0$
of $R_0$ such that ${\cal P}(E_0)=[t_-,t_+]$.
\end{lemma}

\noindent
{\bf Proof of Theorem \ref{teo5} 1.} The boundedness of $S_I$ for each bounded
interval $I$ is trivial -- indeed, if we have a sequence of solutions
to the problem which is unbounded in $\linf$, we divide each equation
by the norm of its solution, as we have  already done a number of times, and  we
find a solution which contradicts Theorem \ref{isolated}. Recall the
regularity result in Theorem \ref{regul}.

For each $n\in\N$ we define $E_n=\calc(\lambda,n)$
as  the connected set given in Lemma \ref{lemche1}.
Then, by Lemma \ref{conn}, there are closed connected subsets $E^N_n$ of
$
\{(u,t)\in E_n\,|\, t\in [-N,N]\}$, for $1\le N\le n$, such that
${\cal P}(E^N_n)=[-N,N]$ and $E^N_n\subset E^{N+1}_n$, for $N=1,2,...,n-1$.
In order to get the last property, we  proceed step by step, defining
$E^N_n$ through Lemma \ref{conn}, by decreasing $N$ starting from $n$.
Now we define the sets $E^N$, for $N\in \N$, as follows :
\begin{eqnarray*}
E^N&=&\{ (u,t)\in C(\overline\Omega)\times\R\,|\,\, \mbox{there exist}\,\,
(u_{\ell_k},
t_{\ell_k})\in E^N_{\ell_k},\\
& & \quad {\ell_k\ge k},\, \forall k\in\N,  (u_{\ell_k},
t_{\ell_k})\to (u,t), \mbox{ as } k\to\infty
\}.
\end{eqnarray*}
We notice that $E^N$ is closed and ${\cal P}(E^N)=[-N,N]$. Since the pairs
$(u,t)\in E^N_n$ are solutions of
$$
F(D^2u,Du,u,x)+\lambda u= t \varphi_1^++h
\quad \mbox{in}\quad \Omega,\quad  u=0 \quad\mbox{on} \quad \partial \Omega,\quad t\in[-N,N],
$$
we see that the set $E^N$ is comprised of solutions of these equations,
and consequently it is compact in $C^1(\overline{\Omega})$.
Then it is easy to see
that for all $\varepsilon>0$ there exists $n_0 \in\N$ such that
$E^N_n\subset B(E^N, \varepsilon)$ for all $n\ge n_0$. Here we denote by
$B(U, \varepsilon)$ the $\varepsilon$-neighborhood of the set $U$.
Indeed, if  there exists
$\varepsilon>0$ and a sequence $\ell_k\ge k$, such that
$(u_{\ell_k},t_{\ell_k})\in
E^N_{\ell_k}\setminus B(E^N, \varepsilon)$, then $t_{\ell_k}$ and
$u_{\ell_k}$ are bounded, and  a subsequence of
$(u_{\ell_k},t_{\ell_k})$ converges to some $(u,t)$ in $E^N$, which is a contradiction.

By the convergence property just proved, we see that $E^N$ is
connected. In fact, if it is not connected, there exist non-empty
closed subsets $U,V$
of
$E^N$ such that $U\cap V=\emptyset$ and $U\cup V=E^N$.
By compactness, there exists $\varepsilon>0$ such that dist$(U,V)>\varepsilon$,
and then $B(U,\varepsilon/4)\cap B(V,\varepsilon/4)=\emptyset$ which
is impossible, since the connected set $E^N_n$ is contained in
$B(U,\varepsilon/4)\cup B(V,\varepsilon/4)$ for $n$ large enough, as stated
in the claim above.

We observe that, according to our construction of the sets $E^N_n$ and $E^N$,
we have $E^N\subset E^{N+1}$ for all $N\in\N$.
So, to complete
the proof of Part 1. we just need to define  $\calc=\calc(\lambda)=\cup_{N\in\N} E^N$, which is
clearly a closed connected set of solutions and ${\cal P}(\calc)=\R$.\hfill $\Box$

\medskip

Before proceeding to  the proof of   Part 2. of Theorem \ref{teo5}, we give a generalized version of the Antimaxmum Principle for fully nonlinear equations, recently proved in \cite{A}.

\begin{prop}\label{antimax} Let  $f\in L^p(\Omega)$, $p>N$, be such that  $f\le 0$, $f\not\equiv0$ in $\Omega$.

 1. There is $\varepsilon_0>0$ (depending on $f$) such that any solution of the equation
\begin{equation}\label{equne} F (D^2u,Du,u,x)+\lambda u=kf\quad\mbox{ in }\; \Omega,\qquad
u=0\quad\mbox{ on }\; \partial\Omega, \end{equation}
satisfies $u<0$ in $\Omega$, provided $\lambda\in (\lambda_1^-, \lambda_1^-+\varepsilon_0)$ and $k\in (0,\infty)$.

2. Equation \re{equne} has no solutions if $\lambda=\lambda_1^-$ and $k>0$.
\end{prop}

\noindent {\bf Proof.} We first prove statement 2. Suppose there is a solution $u$ of \re{equne} with $\lambda=\lambda_1^-$ and $k>0$. If there exists $x_0\in \Omega$ such that
$u(x_0)<0$, then by Theorem~\ref{simple}  there exists  $k_0>0$
such that
$u=k_0\varphi_1^-$,  a contradiction with $f\not\equiv0$. Therefore $u \geq 0$ in  $\Omega$ and then,
by the strong maximum principle, $u > 0$ in  $\Omega$. The existence of such a function contradicts
Theorem \ref{simple}.

Let us now prove statement 1. Suppose there are sequences  $k_n>0$,
$\lambda_n>\lambda_1^-$, $\lambda_n\to \lambda_1^-$, and $\tilde{u}_n$
of solutions of \re{equne}
such that $\tilde{u}_n$ is positive or zero somewhere in $\Omega$.
Then $u_n =\tilde{u}_n/k_n$ has the same property and solves \re{equne}
with $k=1$. Suppose first that $u_n$ is bounded,  then a subsequence of
$u_n$ converges uniformly to a solution of \re{equne} with $\lambda=
\lambda_1^-$ and $k=1$, a contradiction with the result we already proved in 2.
  If $u_n$ is unbounded, then a subsequence of $u_n/\norm{u_n}$ converges
in $C^1(\overline\Omega)$ to the negative function $\varphi_1^-$,
a contradiction as well. \hfill $\Box$

\medskip

We now prove that the solutions of our equation are negative for small $t$
for $t$ below a certain value.

\begin{lemma}\label{refr} Given $R>0$ there are numbers $\varepsilon>0$
and $\bar t$ such that
for all $\lambda\in [\lambda_1^-, \lambda_1^-+\varepsilon)$, $t\le \bar t$, and $h$ with
$\|h\|_{\lp}\le R$, if $u$ solves the equation
\begin{equation}\label{equre} F (D^2u,Du,u,x)+\lambda u=t\varphi_1^+ + h\quad\mbox{ in }\; \Omega,\qquad
u=0\quad\mbox{ on }\; \partial\Omega, \end{equation}
then $u<0$ in $\Omega$.\end{lemma}

\noindent
{\bf Proof.}
Assuming the result is not true, there are sequences $\{t_n\}$, $\{u_n\}$,
$\{\lambda_n\}$ and $\{h_n\}$ such that $\lambda_n\ge \lambda_1^-$,
$\lambda_n\to \lambda_1^-$, $t_n\to -\infty$, $\|h_n\|_{L^p}\le R$,
$u_n$ is positive or zero  at a point in $\Omega$ and
$$ F (D^2u_n,Du_n,u_n,x)+\lambda_n u_n=t_n\varphi_1^++h_n
\quad\mbox{ in }\; \Omega,\qquad
u_n=0\quad\mbox{ on }\; \partial\Omega,
$$
for all $n\in\N$. Defining $v_n=-u_n/t_n$, we can easily check that if $\{v_n\}$ is bounded then a subsequence of it converges to a solution of $F(v) + \lambda_1^-v = -\varphi_1^+$ in $\Omega$, which is a contradiction with part
2. of Proposition \ref{antimax},  while if $\{v_n\}$ is unbounded then
a subsequence of $v_n/\norm{v_n}$
converges in $C^1(\overline{\Omega})$ to $\varphi_1^-<0$,  a contradiction,
since
these functions are positive or zero somewhere.
\hfill $\Box$

\medskip

\noindent
{\bf Proof of Theorem \ref{teo5} 2.} It remains to analyze the asymptotic behavior of the set $S$.
Take any $ u_t\in S_t$, $t\in \R$. It is clear that there exist constants
 $C_0, T>0$, depending only  on $F$, $\Omega$ and $h$, such that
$\norm{u_t}\ge C_0|t|$ if $|t|\ge T$. Indeed, assuming that  $\{t/\norm{u_t}\}$ is not bounded one easily gets the contradiction $0=\pm \varphi_1^+$, after dividing the equation by $t$ and passing to the limit.

First, suppose for contradiction that there exists a compact set $K\subset \Omega$ and sequences $t_n\to -\infty$, $u_n\in S_{t_n}$, such that $u_{t_n}(x_n)
\ge -c$, for some constant $c$ and some $x_n\in K$. Note that by the previous
lemma we already know that $u_{t_n}<0$ in $\Omega$, for large $n$.
Thus, setting
$v_n= u_{t_n}/\norm{u_{t_n}}$, we have $\norm{v_n}=1$, $v_n<0$ in $\Omega$,  $v_n(x_n)\to 0$ as $n\to\infty$, and
$$
F[v_n] + \lambda v_n = (t_n/\norm{u_{t_n}}) \varphi_1^+ + h/\norm{u_{t_n}}\quad \mbox{ in }\; \Omega, \qquad v_n=0\quad\mbox{ on }\; \partial \Omega.$$
Now, if $t_n/\norm{u_{t_n}}\to 0$, a subsequence of $v_n$ converges to
a nontrivial solution of $F[v] + \lambda v =0$, which is a contradiction
with  $\lambda\in(\lambda_1^-,\lambda_1^-+\varepsilon)$.  On the contrary,
if $t_n/\norm{u_{t_n}}\not \to 0$, then a subsequence of $v_n$
converges uniformly to a  solution of $F(v) + \lambda v =-k\varphi_1^+$
for some $k>0$. In addition $v(x_0)=0$ for some $x_0\in K$,  which is a
contradiction with the antimaximum principle, Proposition \ref{antimax},
provided $\varepsilon<\varepsilon_0(-\varphi_1^+)$, with $\varepsilon_0$ defined in that proposition.

Second, suppose there is a sequence $t_n\to +\infty$ such that $u_{t_n}\le C$,
 for some constant $C$. Then, as above, either $v_n= u_{t_n}/\norm{u_{t_n}}$
converges to a nontrivial  solution of $F(v) + \lambda v =0$,  a
contradiction with Theorem \ref{isolated}, or $v_n$ converges  to a nonpositive
solution of $F(v) + \lambda v =k\varphi_1>0$, which is then negative by
Hopf's lemma. This is a contradiction again, here with the definition of $\lambda_1^-$ and $\lambda>\lambda^-_1$.
Theorem \ref{teo5} is proved. \hfill $\Box$

%

\setcounter{equation}{0}
\section{Resonance at $\lambda=\lambda_1^-$. Proof of Theorem 1.4}\label{sect5}

In this section we study the behavior of the set of solutions
of our problem in the second resonant case, that is,
when $\lambda= \lambda_1^-$. For this purpose we
consider sequences $\{\lambda_n\}$ with  $\lambda_n\in (\lambda_1^-,\lambda_1^-+\varepsilon)$ (everywhere in this section $\varepsilon=L$ will be the number which appears in Theorem \ref{teo5}, found in the previous section), which
converge to $\lambda_1^-$, and we study the asymptotic behavior of the connected
sets $\calc=\calc({\lambda_n})\subset S({\lambda_n})$,
obtained in Theorem
\ref{teo5}.

We  modify the definition of
condition ${\cal P}(s)$ as follows:

\begin{itemize}
\item[${\cal P}(s)$]: {\it
There exist  sequences
$\{\lambda_n\}$, $\{h_n\}$ and $\{u_n\}$ such that $\lambda_n>\lambda_1^-$
for all~$n$,  $\lim_{n\to\infty}\lambda_n=\lambda_1^-$,
$h_n\to h$ in $L^p(\Omega),$
$$
F(D^2u_n,Du_n,u_n,x)+\lambda_n u_n=s\varphi_1^++h_n\quad\mbox{ in }\;\Omega,\qquad u_n=0\quad\mbox{ on }\;\partial \Omega,
$$
and  $\norm{u_n}$ is unbounded.}
\end{itemize}
Since no confusion arises with the definition given in Section 3, we
keep the same notation. As before, ${\cal P}(s)$ is equivalent to
\begin{itemize}
\item[${\cal P}(s):$]
{\it
There exist  sequences
$\{\lambda_n\}$, $\{h_n\}$ and $\{u_n\}$ such that $\lambda_n>\lambda_1^-$
for all $n$,  $\lim_{n\to\infty}\lambda_n=\lambda_1^-$, $h_n\to h$ in $\lp$, $\{u_n\}$ is a sequence of solutions of
$
F(D^2u_n,Du_n,u_n,x)+\lambda_n u_n=s\varphi_1^++h_n,
$
 such  that $\norm{u_n}\to\infty$, and  $$\frac{u_n}{\|u_n\|}\to \fimin<0\qquad \mbox{in }\; C^1(\overline{\Omega}).$$}
\end{itemize}
\medskip

Then  we define, as before,
\begin{equation}
t^*_-=\sup\{t\in \R\,|\, {\cal P}(s) \mbox{ for all } s<t\}.
\end{equation}
The following lemmas are necessary to give sense to this definition.
\begin{lema}\label{t1}
${\cal P}(\bar t)$ implies ${\cal P}( t)$ for all $t<\bar t$.
\end{lema}
\noindent
{\bf Proof.}
Assume that there exists $t_0<\bar t$ such that ${\cal P}(t_0)$ is false.
Since ${\cal P}(\bar t)$ holds, there exist sequences $\{\lambda_n\}$, $\{h_n\}$ and
$\{v_n\}$ such that $\lambda_n>\lambda_1^-$
for all~$n$,  $\lim_{n\to\infty}\lambda_n=\lambda_1^-$, $h_n\to h$ in $L^p(\Omega)$,  the solutions of
$$
F(D^2v_n,Dv_n,v_n,x)+\lambda_n v_n=\bar t \varphi_1^++h_n
\quad \mbox{in}\quad \Omega,\quad  v_n=0 \quad\mbox{on} \quad \partial \Omega,
$$
satisfy $\lim_{n\to\infty} \|v_n\|=\infty$, and
$v_n/\|v_n\|$ converges to $\varphi_1^-<0$ in $C^1(\overline \Omega)$, in other words $v_n\le k_n \fimin$, for some sequence $k_n\to\infty$.
On
the other hand, let $\{u_n\}$ be any sequence such that
$$
F(D^2u_n,Du_n,u_n,x)+\lambda_n u_n= t_0 \varphi_1^++h_n
\quad \mbox{in}\quad \Omega,\quad  u_n=0 \quad\mbox{on} \quad \partial \Omega.
$$
Such a sequence exists thanks to Theorem \ref{teo5}. Since we are assuming that
${\cal P}(t_0)$ is false, $\{\|u_n\|\}$ is bounded, so a subsequence of $\{u_n\}$ converges in $C^1(\overline \Omega)$.

Then $|u_n|\le C|\fimin|$ in $\Omega$, by the boundary Lipschitz estimates (recall $\fimin$ has non-zero normal derivative on the boundary, by Hopf's lemma),  so  the above convergence properties of $v_n$ imply
that for $n$ large  $\psi_n=v_n-u_n<0$
in~$\Omega$.
However, by  (H3) we have $F[\psi_n]\ge F[v_n]-F[u_n]$, so
$$
F(D^2\psi_n,D\psi_n,\psi_n,x)+\lambda_n \psi_n\ge  (\bar t-t_0) \varphi_1^+>0
\quad \mbox{in}\quad \Omega,\quad \psi_n =0 \quad\mbox{on} \quad \partial \Omega,
$$
for large $n$, contradicting the definition of $\lambda_1^-$, since $\lambda_n>\lmin$.\hfill $\Box$

\medskip

Now we prove that $t_-^*$  is a real number. We set $T=\{t\in \R\,|\, {\cal P}(t)\}$.
\begin{lema} The set $T$ is not
empty.
\end{lema}
\noindent
{\bf Proof.}
  Assuming the contrary, we find a sequence $\{t_n\}$  such
that ${\cal P}(t_n)$ is false and $t_n\to-\infty$, which implies the existence
of a sequence  $u_n$ satisfying
$$
F(D^2u_n, Du_n,u_n,x)+\lambda_1^- u_n=t_n\varphi_1^++h\quad \mbox{in}\quad \Omega,
\qquad u_n=0 \quad\mbox{on}
\quad \partial \Omega.
$$
This statement follows from Theorem \ref{conv}, through exactly the same argument as the one used in the proof of Lemma \ref{lemr}. Next we see that $v_n=
-{u_n}/t_n$ is unbounded, since the contrary implies the
existence of a solution to
$
F(D^2v, Dv,v,x)+\lambda_1^- v=-\varphi_1^+$ {in} $\Omega$,
$v=0 $ {on}
$\partial \Omega,
$
which was shown to be impossible in Proposition \ref{antimax}.
Then a subsequence of $u_n/\|u_n\|$ converges in $C^1(\overline \Omega)$ to a solution of the equation
$$
F(D^2w, Dw,w,x)+\lambda_1^- w=0\quad \mbox{in}\quad \Omega,
\qquad w=0 \quad\mbox{on}
\quad \partial \Omega,
$$
which implies that $w=\varphi_1^-$. We conclude that $\max_K u_n\to -\infty$
for each compact $K\subset\Omega$.
To complete the proof, let $v$ be the solution of
$$
F(D^2v, Dv,v,x)=-h\quad \mbox{in}\quad \Omega,
\qquad v=0 \quad\mbox{on}
\quad \partial \Omega.
$$
Then, for $n$ large, the function $\psi=u_n+v$
is negative at some point and satisfies
$$
F(D^2\psi, D\psi,\psi,x)+\lambda_1^-\psi \le t_n\varphi_1^++\lambda_1^-v
\quad \mbox{in}\quad \Omega,
\qquad \psi=0 \quad\mbox{on}
\quad \partial \Omega
$$
(we use $F[\psi]\le F[u_n] + F[v]$ which is a consequence of (H0) and (H3)). The quantity $t_n\varphi_1^++\lambda_1^-v$ is strictly negative for large $n$, so by Theorem \ref{simple}  we have $\psi=k\varphi_1^-$,  for some $k>0$, which is a contradiction with the strict inequality $F[\psi]+\lambda_1^-\psi<0$.
Hence $T\not =\emptyset$.\hfill $\Box$

\begin{lema} There exists $\bar t=\bar t(h)\in \mathbb{R}$ such that for any $t\ge \bar {t}$ we can find $C,\delta>0$ such that if $\norm{\tilde h- h}_{\lp}<\delta$, then all solutions to
$$
F(D^2u, Du,u,x)+\lambda v= t\varphi_1^++\tilde h \quad \mbox{in}\quad \Omega,
\qquad u=0 \quad\mbox{on}
\quad \partial \Omega,
$$
with $\lambda\in [\lambda_1^-, \lambda_1^-+\varepsilon)$ satisfy
$\|u\|\le C$. In particular, the set $T$ is bounded above by $\bar t$, that is, $t_-^*$  is finite.
\end{lema}
\noindent
{\bf Proof.}
Assuming the contrary,
we may find
sequences $t_n\to \infty$ as $n\to\infty$, $\lambda_n^{(m)}\in [\lambda_1^-,
\lambda_1^-+\varepsilon)$, $h_n^{(m)}\to h$ in $\lp$ as $m\to\infty$, for each fixed $n$, and  $\{u_n^{(m)}\}$, such that
$$
F(D^2u_n^{(m)},Du_n^{(m)},u_n^{(m)},x)+\lambda_n^{(m)} u_n^{(m)}= t_n\varphi_1^++h_n^{(m)}
\,\,\, \mbox{in}\,\,\, \Omega,\,\,  u_n^{(m)}=0 \,\,\,\mbox{on} \,\,\, \partial \Omega,
$$
and $\{u_n^{(m)}\}$ is unbounded as $ m\to \infty$, for each $n$. Then, as we have done a number of times already, we can divide the last equation by $\|u_n^{(m)}\|$ and use Theorem \ref{conv}, which implies that, up to a subsequence,
  $u_n^{(m)}/\|u_n^{(m)}\|$
converges in $C^1(\overline{\Omega})$, as $m\to \infty$, to a function $\hat u_n\not\equiv0$, which solves
$$
F(D^2\hat u_n,D\hat u_n,\hat u_n,x)+ \lambda_n \hat u_n= 0
\quad \mbox{in}\quad \Omega,\quad  \hat u_n=0 \quad\mbox{on}
\quad \partial \Omega.
$$
This implies that $ \lambda_n=\lambda_1^-$ and $\hat u_n=\varphi_1^-<0$. Hence $u_n^{(m)}\le k_n^{(m)}\fimin$, for some sequence $\{k_n^{(m)}\}$ such that $k_n^{(m)} \to \infty$ as $m\to\infty$.

 Next, we remark that we can find (thanks to Theorems \ref{abp} and \ref{regul}) a constant $C_0=C_0(h)$  such that for any $g\in L^p(\Omega)$ with $\norm{g}_{L^p(\Omega)}\le \! \norm{h}_{L^p(\Omega)} +1$, if $w$ is a solution of
\begin{equation}\label{frety}
F(D^2w,Dw,w,x)=g
\quad \mbox{in}\quad \Omega,\quad  g=0 \quad\mbox{on} \quad \partial \Omega,
\end{equation}
then $\norm{w}_{W^{2,p}(\Omega)}\le C_0$. This of course implies $w\ge \tilde C_0 \fimin$, for some  $\tilde C_0>0$.

Now, for each $n$ we fix $m(n)$ such that $\lambda_n^{(m(n))}<\lambda_1^-+1/n$,
$h_n:=h^{(m(n))}_n$ satisfies $\norm{h_n-h}_{L^p(\Omega)}\le 1/n$, and $u_n:= u^{(m(n))}_n<w$, for each solution $w$ of \re{frety}. So, in particular, $u_n<v_n$, where
 $v_n$ is the solution of
$$
F(D^2v_n,Dv_n,v_n,x)=h_n
\quad \mbox{in}\quad \Omega,\quad  v_n=0 \quad\mbox{on} \quad \partial \Omega.
$$
Then, we choose  $n$ large enough so that
 $t_n \varphi_1^+>\lambda_1^- v_n$, and we see that
 the function $\psi_n=u_n-v_n<0$ satisfies
$$
F(D^2 \psi_n,D  \psi_n, \psi_n,x)+ \lambda_1^- \psi_n\ge t_n \varphi_1^+-\lambda_1^- v_n>0
\quad \mbox{in}\quad \Omega,\quad   \psi_n=0 \quad\mbox{on}
\quad \partial \Omega.
$$
By Theorem \ref{simple} we find that $ \psi_n=\varphi_1^-$, which
 contradicts the last strict inequality.\hfill
$\Box$

\medskip

The next result contains Part 1. in Theorem \ref{teo4}.
\begin{proposition}\label{pp}
The equation
$$
F(D^2u,Du,u,x)+\lambda_1^- u= t \varphi_1^++h
\quad \mbox{in}\quad \Omega,\quad  u=0 \quad\mbox{on} \quad \partial \Omega,
$$

(i) has at least one solution if $t>t^*_-$ ;

(ii) does not have a solution if
$t<t_-^*$.
\end{proposition}

\noindent
{\bf Proof.} (i) is proved in exactly the same way as Proposition \ref{lema33} 1.,
using Theorems \ref{teo5} and \ref{conv}.
In order to prove (ii),
let $t_1\in (t,t^*_-)$. By Lemma \ref{t1}  ${\cal P}(t_1)$ holds, then
there exist sequences
$\{v_n\}$, $\{h_n\}$ and $\{\lambda_n\}$ such that $\{{v_n}\}$ is unbounded,
$\lambda_n> \lambda_1^-$,
$\lambda_n\to \lambda_1^-$, $h_n\to h$, $v_n/\|v_n\|\to \varphi_1^-$ in $C^1(\overline \Omega)$, and
$$F(D^2 v_n, Dv_n,v_n,x)+\lambda_n v_n =t_1\varphi_1^++h_n \quad \mbox{in}\quad
\Omega,\quad  v_n=0 \quad\mbox{on} \quad \partial \Omega. $$

Now, supposing (ii) is false, let $u$ and $w_n$ be  solutions of
 \begin{eqnarray*}
F(D^2 u, Du,u,x)+\lambda_1^- u &=&t\varphi_1^++h \quad \mbox{in}\quad
\Omega,\quad  u=0 \quad\mbox{on} \quad \partial \Omega,\\
F(D^2 w_n, Dw_n,w,x) &=&h_n-h \quad \mbox{in}\quad
\Omega,\quad  w_n=0 \quad\mbox{on} \quad \partial \Omega.
\end{eqnarray*}
Notice that  $w_n\to 0$ in $C^1(\overline{\Omega})$. Then,
for large $n$ we have  $v_n<0$, $\lambda_1^- v_n>\lambda_nv_n$, $v_n-w_n- u<0$ in $\Omega$,   and
\begin{equation}\label{fdsq}
F[v_n-w_n-u]+\lambda_1^-(v_n-w_n-u) \ge
(t_1-t)\varphi_1^+/2>0\quad\mbox{ in }\;\Omega,
\end{equation}
where we used (H3) which implies $F[v_n-w_n-u]\ge F[v_n] -F[w_n] - F[u]$. Then, by Theorem \ref{simple} once more,
we have $v_n-w_n-u=k_n \varphi_1^-$
for some number $k_n\ge0$,
in  contradiction with the strict inequality in \re{fdsq}.
\hfill $ \Box$

\medskip

The next lemma containts statement 3. in Theorem \ref{teo4}.
\begin{lema}\label{tuklem} For each compact interval $I\subset
(t^*_-,\infty)$ there exists a constant $C$, such that for all $\lambda\in [\lambda_1^-, \lambda_1^-+\varepsilon)$ and all $t\in I$, if $u$ is a solution to
$$
F(D^2u, Du,u,x)+\lambda v= t\varphi_1^++h \quad \mbox{in}\quad \Omega,
\qquad u=0 \quad\mbox{on}
\quad \partial \Omega,
$$
then  $\norm{u}_{W^{2,p}(\Omega)}\le C$.
\end{lema}
\noindent
{\bf Proof.}
Recall we already proved in the previous section  that the set of solutions is bounded for $t$ in a bounded interval, provided $\lambda$ is away from the eigenvalue $\lmin$.
Hence if the statement of Lemma \ref{tuklem} is false, then we can
find sequences $t_n\to t_0$ with $t_0>t_-^*$,
$\lambda_n\to \lambda_1^-$ ($\lambda_n= \lambda_1^-$ is allowed), $\{u_n\}$ with $\norm{u_n}\to \infty$
and  $u_n/\norm{u_n}\to \fimin$, such that
$$
F[u_n] + \left( \lambda_n + \frac{1}{\norm{u_n}^2}\right) u_n =
t_0\fipl + h +(t_n-t_0)\fipl + \frac{1}{\norm{u_n}}\left(\frac{u_n}{\norm{u_n}}\right) = t_0\fipl + h_n.
$$
Clearly $h_n\to h$ in $\lp$, so the existence of such a sequence contradicts the definition of the number $t^*_-$ and $t_0>t^*_-$. \hfill $\Box$

\medskip

Before continuing, we set up some notation. The set of solutions
${\cal C}$ found in  Theorem \ref{teo5} will be denoted by
${\cal C(\lambda})$, remembering we work with the equivalent
equation \equ{princ}.
We define the function
${\cal Q}: C(\overline \Omega)\times \R\to\R$, as ${\cal Q}(u,t)=\|u\|$ for $(u,t)\in
C(\overline \Omega)\times \R$, and we recall that
${\cal P}$ is the projection
${\cal P}(u,t)=t$. In the proof of Theorem \ref{teo4} the function
${\cal Q}$ plays a role similar to that of ${\cal P}$ in the proof of Theorem \ref{teo5}.
The following lemma will be needed later.
\begin{lema}\label{NN}
Given $t_1>t^*_-$, there  exists $N_0\in\N$ such that  for every
$\lambda\in (\lambda_1^-,\lambda_1^-+\varepsilon)$ and  $N>N_0$
$$
N\in {\cal Q}(\calc(\lambda)_{[t_1,\infty)})\cap {\cal Q}(\calc(\lambda)_{(-\infty,t_1]}),
$$
that is, for all $\lambda$ larger than and sufficiently close to $\lmin$ and all $N$ large we can find  $u^\prime$, $u^{\prime\prime}$ such that
$$
\|u^\prime\|=\|u^{\prime\prime}\| =N, \quad\mbox{ and}
$$
$$
F[u^\prime] + \lambda u^\prime = t^\prime \fipl + h ,\qquad F[u^{\prime\prime}] + \lambda u^{\prime\prime} = t^{\prime\prime} \fipl +h\quad \mbox{ in }\; \Omega,
$$
 where $t^\prime \ge t_1$ and  $t^{\prime\prime} \le t_1$.
\end{lema}
\noindent
{\bf Proof.} Given $t_1>t^*_-$, we let
$N_0\in\N$ be  an upper bound of the set $\calc(\lambda)_{t_1}$, uniformly in the
interval $\lambda
\in (\lambda_1^-,\lambda_1^-+\varepsilon)$ - such a bound exists by the previous lemma.
The conclusion follows from  Theorem \ref{teo5}, since the set $\calc(\lambda)$
is connected and  the sets  $\calc(\lambda)_t$ contain elements whose norms grow arbitrarily,  as
$t\to\infty$ and as $t\to-\infty.$
\hfill $\Box$

\medskip

\noindent
{\bf Proof of Theorem \ref{teo4}.}
The proof follows an idea similar to the one used in the proof of Theorem \ref{teo5}, but here
we take as a parameter the norm of the solution, instead
of  $t$.

Fix $t_1>t_-^*$. We start with a sequence $\{\lambda_n\}$ with $\lambda_n>\lambda_1^-$
and $\lambda_n\to \lambda_1^-$ as $n\to \infty$.
Then we look at the  connected set of solutions $\calc({\lambda_n})$ given by
Theorem \ref{teo5}, and we take $N\in \N, N> N_0$, where $N_0$ is the number form Lemma \ref{NN}.

By an argument
similar to the one given in the previous section (using Lemma \ref{conn} and Lemma \ref{NN}),
we find that for each
$N=n, n-1,...,N_0+1, N_0$, there is a closed connected subset $E_{n}^N\subset
\{ (u,t)\in \calc(\lambda_n)\,|\, \|u\|\le N\}$ such
that
$$
{\cal Q}((E_{n}^N)_{[t_1,\infty)})=[N_0,N]\quad\mbox{and}\quad
{\cal Q}((E_{n}^N)_{(-\infty,t_1]})=[N_0,N],
$$ for
$N=n,n-1,...,N_0 $. For each $n\in\N$ we construct the sets $E_{n}^N$,
starting with $N=n$ and successively going  down to $N=N_0$.
Thus  $E_{n}^N
\subset E_{n}^{N+1}$, $N= n-1,...,N_0+1, N_0$.
Then we define
\begin{eqnarray*}
E^N&=&\{ (u,t)\in C(\overline\Omega)\times\R\,|\,\, \mbox{there exists}\,\,
(u_{\ell_k},
t_{\ell_k})\in E^N_{\ell_k},\\
& & \quad {\ell_k\ge k},\, \forall k\in\N,  (u_{\ell_k},
t_{\ell_k})\to (u,t), \mbox{ as } k\to\infty
\}.
\end{eqnarray*}
We notice that $E^N$ is closed,
$${\cal Q}((E^N)_{[t_1,\infty)})=[N_0, N]\quad\mbox{and}\quad
{\cal Q}((E^N)_{(-\infty,t_1]})=[N_0, N].
$$
 Since the pairs
$(u,t)\in E^N_{n}$ are solutions of
$$
F(D^2u,Du,u,x)+\lambda_n u= t \varphi_1^++h
\quad \mbox{in}\quad \Omega,\quad  u=0 \quad\mbox{on} \quad \partial \Omega,
$$
the bounded in $\linf$ set $E^N$ is made of solutions of such an equation, but with
$\lambda_1^-$ instead of $\lambda_n$, and consequently $E^N$ is compact.
By a similar argument as the one in the proof of  Theorem \ref{teo5}, we can prove that $E^N$ is
connected. Since, according to our construction, we have  that
$E^N_{n}\subset E^{N+1}_{n}
$ for all $n$, we see that
$E^N\subset E^{N+1}$ for all $N\in\N$. Thus the set
  $\calc=\cup_{N\in\N} E^N$ is a closed connected set of solutions
and
\begin{equation}\label{derf}
{\cal Q}(\calc_{[t_1,\infty)})
=[N_0,\infty)\quad \mbox{and}\quad {\cal Q}(\calc_{[-\infty,t_1]})
=[N_0,\infty).
\end{equation}
Next
 we observe that by the
definition of $t^*_-$ and \re{derf}  we have
${\cal P}(\calc_{[t_1,\infty)})=[t_1,\infty)$.
On the other hand, by Proposition \ref{pp} we know that
$$
(t^*_-,t_1]\subset {\cal P}(\calc_{[-\infty,t_1]})
\subset [t^*_-,t_1],
$$
so that we also have ${\cal Q}(\calc_{[t_-^*,t_1]})
=[N_0,\infty)$. This completes the proof of statement 2. and the first statement in 5. of Theorem \ref{teo4}.

Let us  look at the asymptotic behavior of the set of solutions $S$, as $t\to
\infty$. First,
it is easily proved that if $(u_t,t)\in S$ then $\lim_{t\to\infty}\|u\|=\infty$ (if not, we divide the equation by $t$ and pass to the limit $t\to\infty$, as before).
Suppose now that there is a sequence $t_n\to +\infty$ such that for some
$u_{t_n}\in S_{t_n}$ we have  $u_{t_n}\le C$ for some constant $C$.
Then, as in the proof of Theorem \ref{teo5} 2.,
either  $v_n:= u_{t_n}/\norm{u_{t_n}}$ converges  to a non-positive
solution $v$ of $F[v] + \lambda_1^- v =k\varphi_1^+>0,$ with
$v=0$ on $\partial \Omega$, which is
negative by Hopf's lemma, providing
 a contradiction with Theorem \ref{simple}, or
 $v_n$ converges to a nontrivial  solution of
$F[v] + \lambda_1^- v =0$, with
$v=0$ on $\partial \Omega$.
In this case
%
 $v_n$ converges  to $\fimin<0$ in $C^1(\overline\Omega)$, which implies that
for some sequence $k_n\to \infty$ we have $u_n\le -k_n \varphi_1^+$ in
$\Omega$. Let now $w$ be the solution of $F[w] = h$ in $\Omega$, $w=0$ on
$\partial \Omega$. Then by $\|w\|\le C$ and the Lipschitz estimates  we have $u_n-w<0$ in $\Omega$ if $n$ is sufficiently large, so
\begin{eqnarray*} F[u_n-w] + \lambda_1^- (u_n-w) &\ge& F[u_n]+\lambda_1^- u_n
-(F[w] + \lambda_1^- w)\\
&\ge& t_n\varphi_1^+ - \lambda_1^- w.\end{eqnarray*}
However, the last quantity is positive if $n$ is sufficiently large, yielding
 a contradiction with Theorem \ref{simple}. This gives statement 4 in Theorem \ref{teo4}.

Next we see that there is $R>0$ so that if $(u,t)\in S$, with $t\in [t_-^*,t_1]$ and
$\|u\|_\infty\ge R$, then $u<0$. In fact, if the contrary is true,
then there is a sequence $(u_n,t_n)\in S$, with $t_n\in [t_-^*,t_1]$,
$\|u_n\|_\infty\to\infty$, and such that $u_n$ is positive or zero somewhere in $\Omega$.
But this is impossible since  a subsequence of $u_n/\|u_n\|$ converges in
$C^1(\overline\Omega)$ to $\varphi_1^-$, which is  negative. By the same argument we have $\max_K u_n\to-\infty$ for each $( u_n,t_n)\in S$ such that  $\norm{u_n}\to\infty$ and $t_n\in [t_-^*,t_1]$.
This completes the proof of statement 5 in Theorem \ref{teo4}.

We now turn to the proof of statement 6.
Assume the equation
$$
F(D^2u,Du,u,x)+\lambda_1^- u= t^*_- \varphi_1^++h
\quad \mbox{in}\quad \Omega,\quad  u=0 \quad\mbox{on} \quad \partial \Omega,
$$
has an unbounded set of solutions, that is  $S_{t^*_-}$ is unbounded. Let $u_1, u_2\in
S_{t^*_-}$, then there exists
$R_1>0$ so that whenever $u\in S_{t^*_-}$ and
$\|u\|\ge R_1$ we have $u=u_1+k_1 \varphi_1^-$, for some $k_1>0$. In fact,
we already
know that if $\|u\|$ is large enough
then $u/\|u\|$ is close in $C^1(\overline\Omega)$
to $\varphi_1^-$ and then $\psi=u-u_1<0$ in $\Omega$. Since $\psi$ satisfies
$$
F(D^2\psi,D\psi,\psi,x)+\lambda_1^- \psi\ge 0
\quad \mbox{in}\quad \Omega,\quad  \psi=0 \quad\mbox{on} \quad \partial \Omega,
$$
Theorem \ref{simple} implies $\psi=k_1\varphi_1^-$. In the same way we get
$u=u_2+k_2 \varphi_1^-$ if $\norm{u}\ge \max\{R_1,R_2\}$, for some $R_2>0$,
so $u_1-u_2 = (k_2-k_1)\fimin$.

Finally we prove that if $u + k_1\fimin$ and $u + k_2\fimin$ are in
$S_{t}$ for some $k_2>k_1>0$, then $u + k\fimin\in S_{t}$ for each $k\in(k_1,k_2)$.
This is a simple consequence of the convexity and the homogeneity of $F$.
Indeed, setting $\tilde F=F+\lambda_1^-$,
\begin{eqnarray*}
t\fipl + h &=& \tilde F[u_* + k_1\fimin] + (k-k_1)\tilde F[\fimin]\ge
\tilde F[u_* +k_1\fimin + (k-k_1)\fimin]\\ &=& \tilde F[u+k\fimin]=\tilde F[u_* +k_2\fimin - (k_2-k)\fimin]\\
&\ge&\tilde  F[u_* +k_2\fimin] -(k_2-k)\tilde F[\fimin] =t\fipl + h.
\end{eqnarray*}
Theorem \ref{teo4} is proved. \hfill $\Box$

\setcounter{equation}{0}
\section{ Proof of Theorem \ref{teo6}}\label{sect6}

The proof of Theorem \ref{teo6} relies on an estimate on the difference between  the first eigenvalue of an operator on a domain and a proper subset of the domain, which was proved in \cite{BNV} (Theorem 2.4 in that paper) in the context of general linear operators. We give here a nonlinear version of this result.

Given a smooth bounded domain $A\subset \Omega$, we write $\lambda_1^+(A)$ for the first eigenvalue of the operator $F$ on $A$.
\begin{prop}\label{teo8}
Assume  (H0)-(H3). Let $\Gamma$ be a closed set in $\Omega$, such that $|\Gamma|\ge \alpha_0>0$. Then there exists a constant $\alpha>0$ depending only on $\lambda,\Lambda,N,\gamma,\delta, \Omega, \alpha_0$, such that for any smooth subdomain $A$ of $\Omega\setminus\Gamma$ we have
$$
\lambda_1^+(A)\ge \lambda_1^+(\Omega)+\alpha.
$$
\end{prop}

 The proof of Proposition \ref{teo8} is very similar
to the proof of Theorem 2.4 in  \cite{BNV}. Below we will mention the
points where some small changes have  to be made, but before doing that we show  how we get the proof of Theorem~\ref{teo6}, assuming Proposition
\ref{teo8}.

\medskip

\noindent{\bf Proof of  Theorem \ref{teo6}}. We take $d_0=\alpha/2$, where $\alpha$ is the number from Proposition~\ref{teo8}, with $\alpha_0=|\Omega|/2$.
Suppose for contradiction that we have two different solutions $u_1$ and $u_2$ of \equ{first}, with $F$ satisfying the hypothesis of Theorem \ref{teo6}. We distinguish two cases.

First, suppose the function $w=u_1-u_2$ has a constant sign in $\Omega$,
 say $w\le 0$ (otherwise we take $w=u_2-u_1$). Then (H3) implies  $F(w) \geq 0$
in $\Omega$ and then  $w<0$ in $\Omega$, by Hopf's Lemma.
 The existence of such a function contradicts the definition of $\lmin(\Omega)$ and the assumption $\lmin(\Omega)<0$, see Theorem~\ref{simple}.

Second, if $w=u_1-u_2$ changes sign in $\Omega$, then the sets
 $\Omega_1=\{x \in \Omega \;|\; u_1(x)>u_2(x) \}$ and
$\Omega_2=\{x \in \Omega \;|\; u_2(x)>u_1(x) \}$ are not empty. One of
these sets, say $\Omega_1$, satisfies $|\Omega_1| \leq |\Omega|/2$.
Take $\tilde{\Omega}_1$ to be any connected component of $\Omega_1$ and $A$ to be any smooth subdomain of $\tilde{\Omega}_1$. Then the choice of $d_0$, Proposition \ref{teo8} and $\lpl(\Omega)\ge -d_0$ imply
$$
\lpl(A)\ge {\alpha}/2>0.
$$
Take a sequence of smooth domains  $A_n\subset \tilde{\Omega}_1$ which converges to $\tilde{\Omega}_1$. Then $\lpl(A_n)\ge{\alpha}/2>0$, so by applying the ABP inequality (Theorem \ref{abp}) to $F(w)\ge 0$ in $A_n$ we get
$$
\sup_{A_n} w \le C \sup_{\partial A_n} w.
$$
Letting $n\to \infty$ implies $w\le 0$ in $\tilde{\Omega}_1$, since $w=0$ on $\partial \tilde{\Omega}_1$. This is a contradiction with the definition of $\Omega_1\not=\emptyset$ and proves Theorem \ref{teo6}. \hfill $\Box$

\medskip

\noindent{\bf Proof of Proposition \ref{teo8}}. We follow the proof of Theorem 2.4 in \cite{BNV}, given in Section 9 of that paper. We write $$F(M,p,u,x) = F(M,p,u,x) -\delta u + \delta u =: F_0(M,p,u,x) + \delta u,$$
so that $F_0$ is a proper operator. The operator $F$ plays the role of $L$ in
 \cite{BNV}, $F_0$ plays the role of $M$, $\delta$ replaces $c$, and
we let $q=1+\delta$, as in \cite{BNV}. As shown in \cite{CCKS},
the ABP inequality holds for $F_0$, with a constant which depends only on
$\lambda,\Lambda,\gamma,\delta$ and diam$(\Omega)$.

In what follows we list the results in \cite{BNV} which lead to Proposition
\ref{teo8} and we only note the changes  needed in order to cover the nonlinear case.

Theorem 9.1 in \cite{BNV} is proved in the same way here, but we have to choose
$\sigma>0$ so that $G(D^2e^{\sigma x_1},De^{\sigma x_1},e^{\sigma x_1},x)\ge 1$ -- recall $G$ is defined in (H3) -- which  is easily seen to be possible,
by  (H1), and then we use the inequality $F(M-N, p-q, u-v,x)\le
F(M,p,u,x) - G(N,q,v,x)$, which follows from hypothesis (H3).

The proof of Lemma 9.1 in \cite{BNV} is identical in our situation, as
is the proof of Lemma 9.2, provided we have the concavity of
$\lpl (F_0 + \delta,\Omega)$ in $\delta$, for any proper operator
$F_0$ satisfying our hypotheses, see below.

Theorems 9.2 and 9.3 from \cite{BNV} are well-known  to hold for strong solutions, which is actually the only case in which we use them, if the operators  in their statements are replaced by the operator
 $${\cal L}[u] = \mm(D^2 u ) -\gamma |u| -\delta |u|,$$
which appears in the left-hand side of (H1) -- simply because  ${\cal L}[u]$ is equal to a linear operator acting on $u$, whose coefficients depend on $u$ but their bounds do not. Extensions of these theorems to viscosity solutions can be found in \cite{Wa}, \cite{BS} and in the appendix of \cite{QS3}.

Corollary 9.1 from \cite{BNV} is proved identically here. Further, we need
to modify the proof of Proposition 9.3 in \cite{BNV} in the following way: we take $\nu$ to be the solution of
$$
G(D^2\nu,D\nu,\nu,x) -q\nu = -\chi_\Gamma \quad\mbox{ in }\; \Omega, \qquad \nu = 0 \quad \mbox{ on }\;\partial \Omega,
$$
where $\Gamma$ is as defined in Proposition 9.3 in \cite{BNV}.
We easily check that $G[\cdot] -q\cdot$ is proper, $G[u] -qu\le G[u]
\le F[u]\le 0$ in $\Omega\setminus \Gamma$,
$$
F[u-t\nu]\le F[u]-tG[\nu]\le -tG[\nu]= -tq\nu\le -t\nu$$
 in $\Omega\setminus \Gamma$, and the rest of the proof of Proposition 9.3 is the same.

Finally, Proposition \ref{teo8} follows from the above in exactly the same way
 as Theorem 2.4 in \cite{BNV} follows from Proposition 9.3 there. \hfill $\Box$

\medskip

For completeness we shall briefly sketch the elementary proof of fact that  $\lpl (F_0 + \delta,\Omega)$ is concave in $\delta$. Note that we can repeat exactly the same reasonings as the ones given on pages 50 and 68 of \cite{BNV}, the only difference being that here we need to have the convexity in $z$ of the operator
$${\cal F}(z)(x) = F_0(D^2z + Dz\otimes Dz, Dz, 1,x).$$
This is the content of the following lemma.

\begin{lema} Suppose $F=F(M,p,u)$ satisfies (H0), (H1) and (H3), and let $l:\mathbb{R}^N\to {\cal M}_N(\mathbb{R})$ be a linear map. Then the function
$$
h(p) := F(l(p) + p\otimes p, p, 1) : \mathbb{R}^N\to \mathbb{R}
$$
is convex. \end{lema}

\noindent{\bf Proof.} Suppose $F$ depends only on $M$. Then (H3) implies $F(M)-F(N_1)-F(N_2)\le F(M-N_1-N_2)$, so for any $t\in[0,1]$ and any $p_1,p_2\in\mathbb{R}^N$
\begin{equation}\label{rrrf}
\begin{array}{l}
h(t p_1 + (1-t)p_2)-th(p_1) - (1-t)h(p_2) \\
\le F\left( (t p_1 + (1-t)p_2)\otimes (t p_1 + (1-t)p_2) -t p_1\otimes p_1 -(1-t) p_2\otimes p_2\right).
\end{array}
\end{equation}
By the ellipticity of $F$ it is enough to show that the argument of $F$ in the last inequality is a semi-negative definite matrix. Since $p\otimes q$ is linear in both $p$ and $q$, this is trivially seen to be equivalent to the semi-positive definiteness of
$$
(t-t^2)(p_1\otimes p_1 + p_2\otimes p_2 -  p_1\otimes p_2- p_2\otimes p_1),
$$
that is, of $(t-t^2)((p_1-p_2)\otimes (p_1 -p_2))$, which  is of course true, since $t\in[0,1]$ and the eigenvalues of $q\otimes q$ are $0,\ldots,0,|q|^2$, for each $q\in \mathbb{R}^N$.

If $F=F(M,p,u)$ we have exactly the same reasoning, since in \re{rrrf} we get $F(\cdot, 0,0)$. \hfill $\Box$
\bigskip

\noindent
{\bf Acknowledgements:}
P.F. was  partially supported by Fondecyt Grant \# 1070314,
FONDAP and BASAL-CMM projects
 and  Ecos-Conicyt project C05E09.
A. Q. was partially supported by Fondecyt Grant \# 1070264 and USM
Grant \#  12.08.26 and Programa Basal, CMM, U. de Chile.
\medskip

We would like to thank the anonymous referees for many suggestions which permitted to improve the paper.

\bigskip\bigskip

\begin{flushleft}
 Patricio FELMER\\ Departamento de
Ingenier\'{\i}a  Matem\'atica\\ and
Centro de Modelamiento Matem\'atico,  UMR2071 CNRS-UChile \\
 Universidad de Chile, Casilla 170 Correo 3\\
  Santiago, Chile.\\
  e-mail : \verb"pfelmer@dim.uchile.cl"
  \bigskip

Alexander QUAAS\\
Departamento de  Matem\'atica,\\
 Universidad T\'ecnica Santa Mar\'{i}a\\
Casilla: V-110, Avda. Espa\~na 1680\\
 Valpara\'{\i}so, Chile.\\
 e-mail : \verb"alexander.quaas@usm.cl"
\bigskip

Boyan SIRAKOV (corresponding author) \\
UFR SEGMI \\Universit\'e de Paris 10,\\
 92001 Nanterre Cedex, France,
\
\smallskip

and
\smallskip

CAMS, EHESS\\
54 bd. Raspail\\
75006 Paris, France\\
e-mail : \verb"sirakov@ehess.fr"

\end{flushleft}

\end{document}